\documentclass[11pt]{amsart}
\usepackage{amsmath,amsthm, amscd, amssymb, amsfonts}

\newcommand{\ku}{{\Bbbk}}
\newcommand{\Z}{{\mathbb Z}}

\newcommand{\C}{{\mathcal C}}
\newcommand{\F}{{\mathbb F}}

\newcommand{\Ss}{{\mathcal S}}
\newcommand{\Sg}{{\mathfrak S}}

\newcommand{\End}{\mbox{\rm End\,}}
\newcommand{\Aut}{\mbox{\rm Aut\,}}
\newcommand\tr{\operatorname{tr}}
\newcommand\tot{\operatorname{tot}}

\newcommand\ad{\operatorname{ad}}
\newcommand\Hom{\operatorname{Hom}}
\newcommand\Opext{\operatorname{Opext}}
\newcommand\Tot{\operatorname{Tot}}
\newcommand\Map{\operatorname{Map}}
\newcommand{\fde}{{\triangleright}}
\newcommand{\fiz}{{\triangleleft}}
\newcommand{\wC}{\widehat{C}}
\newcommand{\la}{\langle}
\newcommand{\ra}{\rangle}

\numberwithin{equation}{section}\theoremstyle{plain}

\newtheorem{theorem}{Theorem}[section]
\newtheorem{lema}[theorem]{Lemma}
\newtheorem{cor}[theorem]{Corollary}

\newtheorem{prop}[theorem]{Proposition}
\newtheorem{claim}{Claim}[section]

\theoremstyle{definition}
\newtheorem{definition}[theorem]{Definition}
\newtheorem{exa}[theorem]{Example}

\theoremstyle{remark}
\newtheorem{obs}[theorem]{Remark}

\newcommand\id{\operatorname{id}}

\newcommand\boxe{\begin{tabular}{|p{0,1cm}|}
\hline \\ \hline \end{tabular}}

\newcommand\boxee{\begin{tabular}{|p{0,15cm}|}
\hline \\ \hline \end{tabular}}
\newcommand\boxd{\begin{tabular}{|p{0,3cm}|}
\hline \\ \hline \end{tabular}}

\def\pf{\begin{proof}}
\def\epf{\end{proof}}

\theoremstyle{remark}

\begin{document}

\renewcommand{\baselinestretch}{1.2}
\renewcommand{\thefootnote}{}
\thispagestyle{empty}
\title{Braided Hopf algebras arising from matched pairs of groups}
\author{ Nicol\'as Andruskiewitsch and Sonia Natale}
\address{Facultad de Matem\'atica, Astronom\'\i a y F\'\i sica
\newline \indent
Universidad Nacional de C\'ordoba
\newline
\indent CIEM -- CONICET
\newline
\indent (5000) Ciudad Universitaria, C\'ordoba, Argentina}
\email{andrus@mate.uncor.edu, \quad \emph{URL:}\/
http://www.mate.uncor.edu/andrus}
\address{D\' epartement de math\' ematiques et applications
\newline \indent
\' Ecole Normale Sup\' erieure
\newline
\indent 45, rue d'Ulm
\newline
\indent 75230 Paris Cedex 05, France}
\email{Sonia.Natale@dma.ens.fr, \quad \emph{URL:}\/
http://www.mate.uncor.edu/natale}
\thanks{This work was partially supported by CONICET,
CONICOR, Fundaci\'on Antorchas  and Secyt (UNC)}
\date{Revised version of November 14, 2002.}

\begin{abstract} Let $\ku$ be a field. Let also
$(F, G)$ be a matched pair of groups. We give necessary and
sufficient conditions on a pair $(\sigma, \tau)$ of 2-cocycles in
order that the crossed product algebra  and the crossed coproduct
coalgebra $\ku^G {\,}^{\tau} \#_{\sigma} \ku F$ combine into a
braided Hopf algebra. We also discuss diagonal realizations of
such braided Hopf algebras in the category of Yetter-Drinfeld
modules over a finite group.
\end{abstract}

\maketitle

\section{Introduction}

Let $\ku$ be a field and let $F$, $G$ be finite groups. Given a
right action $\fde$ of $F$ on the set $G$ and a cocycle $\sigma:
F\times F \to (\ku^G)^{\times}$, one forms the crossed product
$\ku^G \#_{\sigma} \ku F$. Dually, given a left action $\fiz$  of
$G$ on the set $F$ and a cocycle $\tau: G\times G \to
(\ku^F)^{\times}$, one forms the crossed coproduct  $\ku^G
{\,}^{\tau}\# \ku F$. In general, $R = \ku^G {\,}^{\tau}
\#_{\sigma} \ku F$ is not a Hopf algebra with these multiplication
and comultiplication. A necessary condition is that the actions
$\fiz$, $\fde$ define a matched pair, that is, that they arise
from an exact factorization $\Sigma = FG$. Let us assume that this
is the case. Then $R$ is a Hopf algebra if and only if $\sigma$
and $\tau$ satisfy a further requirement, which can be expressed
as saying that the pair $(\sigma, \tau)$ is a 1-cocycle in certain
complex. See for example \cite{Maext}. The original sources of
this construction are \cite{kac}, \cite{tak}, \cite{maj}.

The starting point of this paper is the following observation:
even if $(\sigma, \tau)$ is \emph{not} a 1-cocycle, $R$ might
admit a structure of a braided Hopf algebra in the sense of
Takeuchi \cite{T}. That is, under certain conditions, there exists
an invertible solution of the braid equation $c: R\otimes R \to R
\otimes R$, such that the structure maps of $R$ commute with $c$
and the comultiplication $\Delta$ is an algebra map, with respect
to the multiplication in $R\otimes R$ twisted by $c$. The main
result of this paper states necessary and sufficient conditions on
the pair $(\sigma, \tau)$ in order that $R$ is a braided Hopf
algebra with respect to a uniquely determined braiding $c$; {\it
cf.} Theorem \ref{cond-sigma-tau}. It turns out that $c$ is
diagonal in the canonical basis of $R$. Furthermore, $R$ is an
extension of $\ku^G$ by $\ku F$ in this case; and the construction
by Andruskiewitsch and Sommerh\"auser explained in \cite[Ch. 3]{S}
is a particular case of the present one.

A pair $(\sigma, \tau)$, making $R$ into a braided Hopf algebra,
will be called a \emph{braided compatible datum} for the matched
pair $\fde: G \times F \to F$, $\fiz: G \times F \to G$. The main
classification result in \cite{S} states that, provided $\ku$ is
algebraically closed of characteristic zero, every cocommutative
cosemisimple braided Hopf algebra over a cyclic group of order $p$
fits into this construction,  for an appropriate matched pair with
$\fde$ trivial, and a braided compatible datum $\sigma$ and $\tau
= 1$.

The main application we have in view  is the construction of new
examples of Hopf algebras via bosonization or Radford biproducts.
Assume that $R$ is a braided Hopf algebra. By general reasons,
there exists always a Hopf algebra $H$ such that $R$, with this
braiding, is a braided Hopf algebra in the category of
Yetter-Drinfeld modules over $H$. We shall say that $R$ is
\emph{realizable} over $H$. It is interesting to determine all
possible Hopf algebras $H$ such that $R$ is realizable over $H$,
for a fixed $R$.

Note that, if the characteristic of $\ku$ does not divide the
orders of $F$ and $G$, then $R$ is semisimple and cosemisimple.
Thus, if this holds and if $H$ is semisimple and cosemisimple, so
is the biproduct $R \# H$.

There are examples of a braided Hopf algebra $R$ and a Hopf
algebra $H$, with $R$ realizable over $H$, but where neither
$\ku^G$ nor $\ku F$ are realizable over $H$. We shall say that the
extension $\ku^G \hookrightarrow R \to \ku F$ is \emph{realizable}
if $R$, $\ku^G$ and $\ku F$ are realizable, and the inclusion and
projection maps in the extension are morphisms in the category as
well.

We study a distinguished class of realizations of the extension
$\ku^G \hookrightarrow R \to \ku F$ over the group algebra $H =
\ku C$ of a finite group $C$; these are the realizations where
both the action and the coaction are diagonal in the canonical
basis of $R$. Our main result in this direction appears in Theorem
\ref{braid-c-chi}; it allows to avoid the lengthy conditions in
Theorem \ref{cond-sigma-tau}. The biproduct $R \# \ku C$ can be
obtained by iterated extensions from group algebras and dual group
algebras.

We present explicit examples of the general construction over the
field $\mathbb C$ of complex numbers.   Let $p$ and $q$ be
distinct prime numbers. In Proposition \ref{trivial-actions} we
give  examples of braided compatible data, in the case where both
actions $\fiz$ and $\fde$ are trivial: we obtain noncommutative
and noncocommutative braided Hopf algebras of dimension $p^4$.
Another family of examples, together with a diagonal realization
over the group $\mathbb Z_p \oplus \mathbb Z_p$, is constructed in
Proposition \ref{eg-p4q}, as a generalization of the examples in
\cite{S}; these are in general not commutative and not
cocommutative of dimension $p^2q$.

The  paper is organized as follows.  In Section \ref{construction}
we present the construction of a braided Hopf algebra as a
bicrossed product. We give a cohomological interpretation of the
required conditions and consider the problem of equivalences of
braided extensions. In Section \ref{diag-real} we look at those
braided Hopf algebras arising from our construction which admit a
diagonal realization over a finite group $C$; explicit examples of
this situation are constructed in Section \ref{eg-diag-real}. 

By suggestion of the referee, we have also included at the end of the paper an appendix where the main constructions are presented in an alternative language.
We point out that this formalism
can be interpreted in  the language of double categories
as defined by Ehresmann.

\subsection{Notation}

All groups are denoted multiplicatively, unless explicitly stated.
If $G$ is a finite group, we denote by $\ku^G$ the algebra of
functions on $G$; and by $\delta_g$ the canonical idempotent
$\delta_g(h) = \delta_{g, h}$, $h\in G$. These form a basis of
$\ku^G$, the dual basis being the basis $(g)_{g\in G}$ of the
group algebra $\ku G$. The center of $G$ is denoted by $Z(G)$ and
the group of homomorphisms $G \to \ku^{\times}$ is denoted by
$\widehat G$.

\section*{Acknowledgement}
This research was done while the authors stayed at the Department
of Mathematics of the \' Ecole Normale Sup\' erieure, Paris.  They
thank Marc Rosso for the kind  hospitality.

\section{Extensions of braided Hopf algebras arising from matched
pairs.}\label{construction}
\subsection{Matched pairs}
We briefly recall the definition of matched pair of groups. See
\cite{Maext} for further details.

Let $F$ and $G$ be finite groups together with a right
action of $F$ on the set $G$,  and
a left action of $G$ on the set $F$
\begin{equation*}
\fiz: G\times F \to G, \qquad \fde: G\times F \to F.
\end{equation*}
We shall assume that these actions satisfy the following conditions:
\begin{align}\label{comp1}
s \fde xy & = (s \fde x) ((s \fiz x) \fde y), \\
\label{comp2} st \fiz x & = (s \fiz (t \fde x)) (t \fiz x),
\end{align}
for all $s, t \in G$, $x, y \in F$. It follows that $s \fde 1 = 1$ and
$1 \fiz x = 1$, for all $s \in G$, $x \in F$.

Such a data of groups and compatible actions is called a {\it matched
pair} of groups. Given finite groups $F$ and $G$, providing them with a pair of
compatible actions is equivalent to finding a group $\Sigma$ together with an
exact factorization $\Sigma = F G$.

We fix from now on a matched pair of groups $\fiz: G\times F \to G$,
$\fde: G\times F \to F$.
We note the following consequence of the compatibility conditions
\eqref{comp1} and \eqref{comp2}, whose proof is straightforward.

\begin{lema}\label{cons} We have, for all $t \in G$ and $y \in F$,

(i) $(t \fiz y)^{-1} = t^{-1} \fiz (t \fde y)$;

(ii) $(t \fde y)^{-1} = (t \fiz y) \fde y^{-1}$. \qed \end{lema}

\subsection{Bicrossed products.}\label{bicpro}

We consider the associated left action of $F$ on $\ku^G$, $(x.
\phi)(g) = \phi(g \fiz x)$, $\phi \in \ku^G$; in particular, $x .
\delta_g = \delta_{g \fiz x^{-1}}$. Let $\sigma: F \times F \to
(\ku^{\times})^G$ be a normalized 2-cocycle. If we write $\sigma =
\sum_{g\in G} \sigma_g \delta_g$, then the cocycle and the
normalizing conditions read, respectively, as follows:
\begin{align}\label{cociclo-sigma}
& \sigma_{g \fiz x} (y,z) \sigma_{g} (x,yz) =  \sigma_{g} (xy,z)
\sigma_{g} (x, y), \\
\label{norm-sigma} & \sigma_g (x, 1) = 1 = \sigma_g (1, x), \qquad g\in
G, x,y,z \in F.
\end{align}

We also consider the associated right action of $G$ on $\ku^F$,
$(\psi .g) (x) = \psi(x \fde g)$, $\psi \in \ku^F$. Let $\tau =
\sum_{x\in F} \tau_x \delta_x: G \times G \to (\ku^{\times})^F$ be
a normalized 2-cocycle; so that we have
\begin{align}\label{cociclo-tau}
& \tau_{x} (gh, k) \tau_{k \fde x} (g,h) =  \tau_{x} (h,k)
\tau_{x} (g, hk), \\
\label{norm-tau} & \tau_x (g, 1) = 1 = \tau_x (1, g), \qquad g,h,k\in
G, x \in F. \end{align}

We endow the vector space $R = \ku^G \otimes \ku F$ with the
crossed product algebra structure and the crossed coproduct
coalgebra structure. By abuse of terminology, we refer to $R$ as a
\emph{bicrossed product}.

We shall use the notation $\delta_g x$ to indicate the element
$\delta_g \otimes x \in R$. Then the multiplication of $R$  is
determined by
\begin{equation}\label{producto-def}
(\delta_g  x)(\delta_h  y) =
\delta_{g \fiz x, h}\, \sigma_g(x, y) \delta_g   xy,
\qquad g,h\in G, x, y\in F;
\end{equation}
and the comultiplication is determined by
\begin{equation}\label{coproducto-def}
 \Delta(\delta_g x) = \sum_{t\in G} \tau_x(t, t^{-1}g )\,
 \delta_t (t^{-1}g \fde x) \otimes \delta_{t^{-1} g} x,
 \qquad g\in G, x \in F.
\end{equation}

In the following lemma we give the necessary normalization conditions
on $\sigma$ and $\tau$ in order that the unit and counit maps preserve
the coalgebra and algebra structures, respectively.

\begin{lema}\label{norm2} (i) $\epsilon \otimes \epsilon: R \to \ku$ is
an
algebra map if and only if
\begin{equation}\label{norm2-sigma} \sigma_1(g, h) = 1,
\qquad \forall g, h \in G; \end{equation}

(ii) $\Delta (1) = 1 \otimes 1$ if and only if
\begin{equation}\label{norm2-tau} \tau_1(x, y) = 1,
\qquad \forall x, y \in F. \end{equation}
\end{lema}

\pf Straightforward.  \epf

We next show that the formula for the antipode still provides the
inverse of the identity, even if $\sigma$ and $\tau$ do not
satisfy any compatibility condition.

\begin{lema} The map $\Ss$ defined by
\begin{equation}\label{antipode}
\Ss(\delta_g  x)=  \sigma_{(g \fiz x)^{-1}}((g \fde x)^{-1}, g \fde
x)^{-1}
\, \tau_x(g^{-1}, g)^{-1}\, \delta_{(g \fiz x)^{-1}}\, (g \fde x)^{-1}
\qquad g\in G, x\in F.
\end{equation}
is the inverse of the identity map with respect to the convolution
product in $\End R$. \end{lema}

\pf Letting $x = z = y^{-1}$ in the cocycle condition
\eqref{cociclo-sigma}, we get
\begin{equation}
\label{comp0} \sigma_{g \fiz x} (x^{-1},x) = \sigma_{g} (x,x^{-1}),
\qquad g\in G, x\in F.
\end{equation}
Combining \eqref{comp0} with Lemma \ref{cons} (i), we have
\begin{equation}
\label{comp0-1} \sigma_{(t \fiz x)^{-1}} ((t \fde x)^{-1},t \fde x) =
\sigma_{t^{-1}} (t \fde x,(t \fde x)^{-1}), \qquad t\in G, x\in F.
\end{equation}
Let $X = \delta_gx \in R$. We compute
\begin{align*}
X_1 \Ss(X_2) & = \sum_{t \in G} \tau_x(gt^{-1}, t) \delta_{gt^{-1}} (t
\fde x) \, \Ss(\delta_t x) \\
& =  \sum_{t \in G} \tau_x(gt^{-1}, t) \, \sigma_{(t \fiz x)^{-1}}((t
\fde x)^{-1}, t \fde x)^{-1}
\, \tau_x(t^{-1}, t)^{-1}\, \delta_{gt^{-1}} (t \fde x) \, \delta_{(t
\fiz x)^{-1}}\, (t \fde x)^{-1} \\
& = \sum_{t \in G} \tau_x(gt^{-1}, t)\,  \sigma_{(t \fiz x)^{-1}}((t
\fde x)^{-1}, t \fde x)^{-1}
\, \tau_x(t^{-1}, t)^{-1}\, \sigma_{gt^{-1}} (t \fde x, (t \fde
x)^{-1}) \,
\delta_{gt^{-1}, (t \fiz x)^{-1} \fiz (t \fde x)^{-1}} \,
\delta_{gt^{-1}} \\
& = \delta_{g, 1} \sum_{t \in G} \delta_t = \delta_{g, 1} 1 = \epsilon
(X) 1.
\end{align*}
In the fourth equality we have used \eqref{comp0-1} and the fact that
$\delta_{gt^{-1}, (t \fiz x)^{-1} \fiz (t \fde x)^{-1}} =
\delta_{gt^{-1}, t^{-1}}$, which follows from Lemma \ref{cons} (i).
This proves the lemma. \epf

\subsection{Braiding}\label{braiding}

We shall consider in this subsection $2$-cocycles $\sigma: F
\times F \to (\ku^{\times})^G$ and $\tau: G \times G \to
(\ku^{\times})^F$ satisfying the  normalization conditions
\eqref{norm-sigma}, \eqref{norm-tau}, \eqref{norm2-sigma} and
\eqref{norm2-tau}.

\bigbreak We shall give necessary and sufficient conditions in
order that the resulting algebra and coalgebra structures on $R$
associated to the data $\fiz$, $\fde$, $\sigma$ and $\tau$ make it
a braided Hopf algebra in the sense of \cite[Section 5]{T}. That
is, we shall determine when there exists an invertible linear map
$c: R \otimes R \to R \otimes R$, which satisfies the following
conditions:

\bigbreak
(i) $c$ is a solution of the braid equation $(c \otimes
\id) (\id \otimes c) (c \otimes \id) = (\id \otimes c) (c \otimes
\id) (\id \otimes c)$;

\bigbreak (ii) the structure maps of $R$ commute with the
braidings. According to the definition in \cite{T}, $m$ commutes
with $c$ if and only if
\begin{align}\label{tak-m1}
c (\id \otimes m) &= (m \otimes \id) c_{1, 2}
\\\label{tak-m2}
(\id \otimes m)c_{2, 1} &= c (m \otimes \id), \end{align} where
the braidings $c_{1, 2}: R \otimes R^{\otimes 2} \to R^{\otimes 2}
\otimes R$ and $c_{2, 1}: R^{\otimes 2} \otimes R \to R \otimes
R^{\otimes 2}$, are defined by $$c_{1, 2} : = (\id \otimes c) (c
\otimes \id ), \qquad c_{2, 1}: = (c \otimes \id ) (\id \otimes
c).$$ Similarly, $\Delta$ commutes with $c$ if and only if
\begin{align}
\label{tak-d1}(\Delta \otimes \id) c &= c_{1, 2} (\id \otimes
\Delta) \\\label{tak-d2}(\id \otimes \Delta) c &= c_{2, 1} (\Delta
\otimes \id).\end{align}

\bigbreak
(iii) $\epsilon: R \to \ku$ is an algebra map,
$\Delta(1) = 1 \otimes 1$, and  $\Delta : R \to R
\underline{\otimes} R$ is an algebra map; here the product in $R
\underline{\otimes} R$ is "twisted" by $c$: $m_{R
\underline{\otimes} R} = (m_R \otimes m_R) (\id \otimes c \otimes
\id )$. Moreover, the identity map has a convolution inverse
$\Ss$, called the antipode.

\bigbreak See \cite[Definition 5.1]{T}. In the case where only
conditions (i) and (iii) are satisfied, we shall say that $R$ is a {\it
prebraided} Hopf algebra.

\begin{obs} In order that the product in $R \underline{\otimes} R$
be associative, we must require that the multiplication map $R
\otimes R \to R$ commutes with the braiding.
\end{obs}

\begin{definition} Let $S$ be a braided Hopf algebra and let $H$ be a
Hopf algebra.
We shall say that $S$ is {\it realizable} over $H$ if it can be
endowed with a left action $H \otimes R \to R$ and a left coaction
$R \to H \otimes R$, such that $S$ is a braided Hopf algebra in
the category ${}^H_H\mathcal{YD}$ of Yetter-Drinfeld modules over
$H$, with the braiding $c$ being the corresponding braiding in
${}^H_H\mathcal{YD}$.
\end{definition}

One may also consider the related notion of being realizable over
a quasitriangular Hopf algebra; this point of view will not be
discussed in this paper.

\bigbreak
Let $S$ be a finite dimensional braided Hopf algebra.
Recall from \cite[Theorem 5.7]{T}, that $S$ is realizable over  a
(non unique) Hopf algebra $H$ if and only if  $S$ is a {\it rigid}
braided Hopf algebra, which means  that the braiding $c$ is rigid.

\bigbreak
 We define $c: R \otimes R \to R \otimes R$ in the form
\begin{equation}\label{def-c}
c(\delta_gx \otimes \delta_hy) =  Q^{x, y}_{g, h}\, \delta_hy
\otimes \delta_gx, \qquad g, h \in G, \, x, y \in F,
\end{equation}
where $Q: G^2 \times F^2 \to \ku^{\times}$ is a map. Note that $c$
is diagonal and thus automatically satisfies the braid equation.
Moreover, since the scalars $Q^{x, y}_{g, h}$ are non-zero by
assumption, $c$ is rigid.

\bigbreak
The following proposition generalizes \cite[Proposition 4.7]{Maext}.

\begin{prop}\label{cond-q}
$(R, c)$ is a prebraided  Hopf algebra if and only if the
following compatibility condition holds, for all $s, t \in G$, $x,
y \in F$:
\begin{equation}\label{comp3}
\sigma_{ts}(x, y) \tau_{xy}(t, s)  =  Q^{x, (s \fiz x) \fde y}_{s,
t \fiz (s \fde x)} \,
 \tau_x(t, s) \, \tau_y(t \fiz (s \fde x), s \fiz x) \,
\sigma_{t}(s \fde x, (s \fiz x) \fde y) \, \sigma_{s}(x, y).
\end{equation} \end{prop}

\pf Note that,  as a consequence of \eqref{norm2-sigma} and
\eqref{norm2-tau}, the compatibility condition \eqref{comp3}
implies the following normalization conditions on $Q$:
\begin{equation}\label{norm-c}
Q^{1, y}_{g, h} =  Q^{x, 1}_{g, h} =  Q^{x, y}_{1, h} = Q^{x,
y}_{g, 1} = 1. \end{equation} We have already established the
existence of an antipode. By Lemma \ref{norm2} the counit is a
morphism of algebras and $\Delta(1) = 1 \otimes 1$. We shall prove
that condition \eqref{comp3} is equivalent to the comultiplication
$\Delta: R \to R \underline{\otimes} R$ being a morphism of
algebras.

 Let $g, h \in G$, $x, y \in F$.
We denote by $\bullet$ the product in $R \underline{\otimes} R$
twisted by $c$. We compute
\begin{align*}
\Delta(\delta_gx) \bullet \Delta(\delta_hy) & = \sum_{s, t \in G}
\tau_x(t, t^{-1}g) \, \tau_y(s, s^{-1}h) \, \left(
\delta_t(t^{-1}g \fde x) \otimes \delta_{t^{-1}g} x \right)
\bullet \left( \delta_s(s^{-1}h \fde y)  \otimes \delta_{s^{-1}h}
y \right) \\ & = \sum_{s, t \in G} \tau_x(t, t^{-1}g) \, \tau_y(s,
s^{-1}h) \, Q^{x, s^{-1}h \fde y}_{t^{-1}g, s} \, \delta_{t, s
\fiz (t^{-1}g \fde x)^{-1}} \, \delta_{t^{-1}g, s^{-1}h \fiz
x^{-1}} \\ & \sigma_{t}(t^{-1}g \fde x, s^{-1}h \fde y) \,
\sigma_{t^{-1}g}(x, y) \, \delta_t (t^{-1}g \fde x) (s^{-1}h \fde
y) \otimes \delta_{t^{-1}g} xy;
\end{align*}
using  the compatibility conditions \eqref{comp1} and
\eqref{comp2}, this equals
\begin{align*}
\delta_{g, h \fiz x^{-1}} & \, \sum_{t \in G} \tau_x(t, t^{-1}g)
\,  \tau_y(t \fiz (t^{-1}g \fde x), (t \fiz (t^{-1}g \fde
x))^{-1}h) \, Q^{x, (t^{-1}h \fiz x) \fde y}_{t^{-1}g, t \fiz
(t^{-1}g \fde x)} \\ & \sigma_{t}(t^{-1}g \fde x, (t^{-1}g \fiz x)
\fde y) \, \sigma_{t^{-1}g}(x, y) \, \delta_t (t^{-1}g \fde x)
((t^{-1}g \fiz x) \fde y) \otimes \delta_{t^{-1}g} xy.
\end{align*}
On the other hand, we have
\begin{equation}\label{delta-m}
\Delta m(\delta_gx \otimes \delta_hy) = \delta_{g, h \fiz x^{-1}}
\sigma_g(x, y) \sum_{s \in G} \tau_{xy}(s, s^{-1}g) \delta_s
(s^{-1}g \fde (xy)) \otimes \delta_{s^{-1}g} xy;
\end{equation}
Hence, we find that $\Delta$ is an algebra map if and only if
\begin{align*}
\delta_{g, h \fiz x^{-1}} \sigma_g(x, y) \tau_{xy}(t, t^{-1}g) & =
\delta_{g, h \fiz x^{-1}} \tau_x(t, t^{-1}g) \, \tau_y(t \fiz
(t^{-1}g \fde x), (t \fiz (t^{-1}g \fde x))^{-1}h) \, Q^{x,
(t^{-1}g \fiz x) \fde y}_{t^{-1}g, t \fiz (t^{-1}g \fde x)} \\ &
\sigma_{t}(t^{-1}g \fde x, (t^{-1}g \fiz x) \fde y) \,
\sigma_{t^{-1}g}(x, y). \end{align*} Letting $s = t^{-1}g$,  this
condition is equivalent to the claimed one. This finishes the
proof of the proposition. \epf

Proposition \ref{cond-q}  allows us to construct, for {\it any}
normalized $2$-cocycles $\sigma$ and $\tau$, a prebraided Hopf
algebra structure on $R$.

\begin{prop}\label{const-struct} There exists a unique braiding
 $c: R \otimes R \to R \otimes R$ making $R$ into a prebraided
 Hopf algebra: it is given by \eqref{def-c}, where $Q: G^2 \times F^2
\to
\ku^{\times}$ is the map defined in the form
\begin{multline}\label{def-q}
Q^{x, y}_{g, h} : = \sigma_{(h \fiz (g \fde x)^{-1})g}\left( x, (g \fiz
x)^{-1} \fde y \right) \,
\sigma_{h \fiz (g \fde x)^{-1}}(g \fde x, y)^{-1} \,
\sigma_g(x, (g \fiz x)^{-1} \fde y)^{-1} \\
\quad \tau_{x ((g \fiz x)^{-1} \fde y)}(h \fiz (g \fde x)^{-1}, g) \,
\tau_{(g \fiz x)^{-1} \fde y} (h, g \fiz x)^{-1} \,
\tau_x(h \fiz (g \fde x)^{-1}, g)^{-1},
\end{multline}
for all $g, h \in G$, $x, y \in F$.

In particular, every braided Hopf algebra structure on $R$ is
realizable over some Hopf algebra $H$.
\end{prop}

An alternative proof of the last statement in the proposition is
given in Lemma \ref{close} below.

\pf It is easy to see that formula \eqref{def-q} is equivalent to
\eqref{comp3}. Therefore, if $c$ is given by \eqref{def-q}, $(R,
c)$ is a prebraided Hopf algebra.

It follows from \cite{Sb}, that the associativity,
coassociativity, unit, counit and antipode axioms on $R$, together
with the condition $\Delta m = (m \otimes m) (\id \otimes c
\otimes \id ) (\Delta \otimes \Delta)$, uniquely determine  the
braiding $c$ by means of the formula $c = (m \otimes m) (\Ss
\otimes \Delta m \otimes \Ss) (\Delta \otimes \Delta)$. Actually,
the argument in \cite{Sb} does not need the associativity of $R
\underline{\otimes} R$.  Therefore, for fixed $\sigma$ and
$\tau$, the braiding $c$ making $R$ into a prebraided Hopf algebra
is unique, and has necessarily the prescribed form.

 In particular, all
such braidings are 'diagonal' in the basis $\delta_gx$, $g \in G$,
$x \in F$, and they are moreover rigid. It follows that every
braided Hopf algebra structure on $R$ is realizable, as claimed.
\epf

\begin{obs} The normalization conditions  \eqref{norm2-sigma},
\eqref{norm2-tau}
in Lemma \ref{norm2} are equivalent to the normalization
conditions on $Q$ in \eqref{norm-c},  in view of \eqref{def-q}.
 \end{obs}

The following lemma gives necessary and sufficient conditions in
order that condition (ii) be satisfied.

\begin{lema}\label{1-cociclos} (i) The multiplication map $m: R \otimes
R  \to R$
commutes with $c$ if and only if
\begin{align}
\label{comm1-m} Q^{x, yz}_{g, s} & = Q^{x, y}_{g, s} \, Q^{x,
z}_{g, s \fiz y}, \\ \label{comm2-m}Q^{xy, z}_{g, s} & =
Q^{x,z}_{g, s} \, Q^{y, z}_{g \fiz x, s}, \qquad \forall x, y, z
\in F, \, g, s \in G.
\end{align}

(ii) The comultiplication map $\Delta: R \to R \otimes R$ commutes
with $c$ if and only if
\begin{align}
\label{comm1-delta} Q^{x, y}_{g, ts} & = Q^{x, s \fde y}_{g, t} \,
Q^{x, y}_{g, s},
\\ \label{comm2-delta} Q^{x, y}_{ts, g} & = Q^{s \fde x, y}_{t, g} \,
Q^{x, y}_{s, g}, \qquad \forall x, y  \in F, \, g, t, s \in G.
\end{align} \end{lema}

\pf (i) Conditions \eqref{tak-m1}, \eqref{tak-m2} are equivalent,
in this case, to \eqref{comm1-m} and \eqref{comm2-m},
respectively.

(ii) Similarly, conditions \eqref{tak-d1}, \eqref{tak-d2}
correspond to \eqref{comm1-delta} and \eqref{comm2-delta},
respectively. \epf

The following theorem is a consequence of Lemma \ref{1-cociclos}
and Proposition \ref{const-struct}. It  gives the necessary and
sufficient conditions on $\sigma$ and $\tau$ in order that $R$ be
a braided Hopf algebra.

\begin{theorem}\label{cond-sigma-tau} Let $c: R \otimes R \to R \otimes
R$
be given by \eqref{def-c}, where $Q: G^2 \times F^2 \to
\ku^{\times}$ is the map defined by \eqref{def-q}. Then $(R, c)$
is a braided Hopf algebra if and only if the following
compatibility conditions hold, for all $g, s, t \in G$, $x, y, z
\in F$:
\begin{flalign*}
(1)\quad  & \sigma_{(s \fiz (g \fde x)^{-1})g}\left( x, (g \fiz
x)^{-1} \fde yz \right) \, \sigma_{s \fiz (g \fde x)^{-1}}(g \fde
x, yz)^{-1} \, \sigma_g(x, (g \fiz x)^{-1} \fde yz)^{-1} & \\ &
\tau_{x ((g \fiz x)^{-1} \fde yz)}(s \fiz (g \fde x)^{-1}, g) \,
\tau_{(g \fiz x)^{-1} \fde yz} (s, g \fiz x)^{-1}  & \\ & =
\sigma_{(s \fiz (g \fde x)^{-1})g}\left( x, (g \fiz x)^{-1} \fde y
\right) \, \sigma_{s \fiz (g \fde x)^{-1}}(g \fde x, y)^{-1} \,
\sigma_g(x, (g \fiz x)^{-1} \fde y)^{-1} & \\ & \tau_{x ((g \fiz
x)^{-1} \fde y)}(s \fiz (g \fde x)^{-1}, g) \, \tau_{(g \fiz
x)^{-1} \fde y} (s, g \fiz x)^{-1}  & \\ & \sigma_{(s \fiz y(g \fde
x)^{-1})g}\left( x, (g \fiz x)^{-1} \fde z \right) \, \sigma_{s
\fiz (g \fde x)^{-1}}(g \fde x, z)^{-1} \, \sigma_g(x, (g \fiz
x)^{-1} \fde z)^{-1} & \\ & \tau_{x ((g \fiz x)^{-1} \fde z)}(s \fiz
y(g \fde x)^{-1}, g) \, \tau_{(g \fiz x)^{-1} \fde z} (s \fiz y, g
\fiz x)^{-1} \, \tau_x(s \fiz y(g \fde x)^{-1}, g)^{-1};
\end{flalign*}

\begin{flalign*}
(2)\quad  & \sigma_{(s \fiz (g \fde xy)^{-1})g}\left( xy, (g \fiz
xy)^{-1} \fde z \right) \, \sigma_{s \fiz (g \fde xy)^{-1}}(g \fde
xy, z)^{-1} \, \sigma_g(xy, (g \fiz xy)^{-1} \fde z)^{-1} & \\ &
\tau_{xy ((g \fiz xy)^{-1} \fde z)}(s \fiz (g \fde xy)^{-1}, g) \,
 \tau_{xy}(s \fiz (g \fde xy)^{-1}, g)^{-1} & \\ & =
\sigma_{(s \fiz (g \fde x)^{-1})g}\left( x, (g \fiz x)^{-1} \fde z
\right) \, \sigma_{s \fiz (g \fde x)^{-1}}(g \fde x, z)^{-1} \,
\sigma_g(x, (g \fiz x)^{-1} \fde z)^{-1} & \\ & \tau_{x ((g \fiz
x)^{-1} \fde z)}(s \fiz (g \fde x)^{-1}, g) \, \tau_{(g \fiz
x)^{-1} \fde z} (s, g \fiz x)^{-1} \, \tau_x(s \fiz (g \fde
x)^{-1}, g)^{-1} & \\ & \sigma_{(s \fiz ((g \fiz x) \fde y)^{-1}) (g
\fiz x)}\left( y, ((g \fiz xy)^{-1} \fde z \right) \, \sigma_{s
\fiz ((g \fiz x) \fde y)^{-1}}((g \fiz x) \fde y, z)^{-1} \,
\sigma_{g \fiz x}(y, (g \fiz xy)^{-1} \fde z)^{-1} & \\ & \tau_{y
((g \fiz xy)^{-1} \fde z)}(s \fiz ((g \fiz x) \fde y)^{-1}, g \fiz
x) \,  \tau_y(s \fiz ((g \fiz x) \fde y)^{-1}, g \fiz x)^{-1};
\end{flalign*}

\begin{flalign*}
(3) \quad & \sigma_{(ts \fiz (g \fde x)^{-1})g}\left( x, (g \fiz
x)^{-1} \fde y \right) \, \sigma_{ts \fiz (g \fde x)^{-1}}(g \fde
x, y)^{-1}  & \\ & \tau_{x ((g \fiz x)^{-1} \fde y)}(ts \fiz (g \fde
x)^{-1}, g) \, \tau_{(g \fiz x)^{-1} \fde y} (ts, g \fiz x)^{-1}
\, \tau_x(ts \fiz (g \fde x)^{-1}, g)^{-1} & \\ & = \sigma_{(t \fiz
(g \fde x)^{-1})g}\left( x, (g \fiz x)^{-1}s \fde y \right) \,
\sigma_{t \fiz (g \fde x)^{-1}}(g \fde x, s \fde y)^{-1}
\sigma_{g}(x, (g \fiz x )^{-1}s \fde y)^{-1}
& \\ & \tau_{x ((g \fiz x)^{-1}s \fde y)}(t \fiz (g \fde x)^{-1}, g)
\, \tau_{(g \fiz x)^{-1}s \fde y} (t, g \fiz x)^{-1} \, \tau_x(t
\fiz (g \fde x)^{-1}, g)^{-1} \\ & \sigma_{(s \fiz (g \fde
x)^{-1})g}\left( x, (g \fiz x)^{-1} \fde y \right) \, \sigma_{s
\fiz (g \fde x)^{-1}}(g \fde x, y)^{-1}
& \\ & \tau_{x ((g \fiz x)^{-1} \fde y)}(s \fiz (g \fde x)^{-1}, g)
\, \tau_{(g \fiz x)^{-1} \fde y} (s, g \fiz x)^{-1} \, \tau_x(s
\fiz (g \fde x)^{-1}, g)^{-1};
\end{flalign*}

\begin{flalign*}
(4) \quad & \sigma_{(s \fiz (gt \fde x)^{-1})gt}\left( x, (gt \fiz
x)^{-1} \fde y \right) \, \sigma_{gt}(x, (gt \fiz x)^{-1} \fde
y)^{-1} & \\ & \tau_{x ((gt \fiz x)^{-1} \fde y)}(s \fiz (gt \fde
x)^{-1}, gt) \, \tau_{(gt \fiz x)^{-1} \fde y} (s, gt \fiz x)^{-1}
\, \tau_x(s \fiz (gt \fde x)^{-1}, gt)^{-1} & \\ & = \sigma_{(s \fiz
(g \fde (t \fde x))^{-1})g}\left( t \fde x, (g \fiz (t \fde
x))^{-1} \fde y \right) \,  \sigma_g(t \fde x, (g \fiz (t \fde
x))^{-1} \fde y)^{-1} & \\ & \tau_{(t \fde x) ((g \fiz x)^{-1} \fde
y)}(s \fiz (gt \fde x)^{-1}, g) \, \tau_{(g \fiz (t \fde x))^{-1}
\fde y} (s, g \fiz (t \fde x))^{-1} \, \tau_{t \fde x}(s \fiz (gt
\fde x)^{-1}, g)^{-1} & \\ & \sigma_{(s \fiz (t \fde
x)^{-1})t}\left( x, (t \fiz x)^{-1} \fde y \right) \, \sigma_{s
\fiz (t \fde x)^{-1}}(t \fde x, y)^{-1} \, \sigma_t(x, (t \fiz
x)^{-1} \fde y)^{-1} & \\ & \tau_{x ((t \fiz x)^{-1} \fde y)}(s \fiz
(t \fde x)^{-1}, t) \, \tau_{(t \fiz x)^{-1} \fde y} (s, t \fiz
x)^{-1} \, \tau_x(s \fiz (t \fde x)^{-1}, t)^{-1}.
\end{flalign*}

This is the case if and only if there exists a Hopf algebra $H$
such that $R$ is a braided Hopf algebra in ${}^H_H\mathcal{YD}$.
\qed \end{theorem}

We stress again that  the Hopf algebra $H$ in Theorem
\ref{cond-sigma-tau} is not unique.

\begin{definition} A pair $(\tau, \sigma)$  of 2-cocycles  $\tau: G
\times G \to
(\ku^{\times})^F$ and $\sigma: F \times F \to (\ku^{\times})^G$,
satisfying conditions \eqref{norm-sigma}, \eqref{norm-tau},
\eqref{norm2-sigma} and \eqref{norm2-tau}, and the compatibility
conditions (1)--(4) in Theorem \ref{cond-sigma-tau} will be called
a \emph{braided compatible datum} for the matched pair $\fde : G
\times F \to$, $\fiz : G \times F \to G$.  \end{definition}

Given a pair $(\sigma, \tau)$, the compatibility conditions (1)--(4) in Theorem
\ref{cond-sigma-tau} are not easy to check. An alternative way to
construct braided Hopf algebras is indicated in Theorem
\ref{braid-c-chi} below.

\begin{obs}  Suppose that $R$ is realizable
over a (finite-dimensional) semisimple Hopf algebra $H$.  Then the
symmetrizations $Q^{x,y}_{g,h}Q^{y, x}_{h, g}$ are roots of unity,
for all $x, y \in F$, $g, h \in G$.

\pf The map $c^2 : R \otimes R \to R \otimes R$ is given by the
action of $\mathcal R_{21}\mathcal R$ on $R \otimes R$, where
$\mathcal R$ is the canonical $R$-matrix of $D(H)$. By
\cite{eg-exp}, the order of $\mathcal R_{21}\mathcal R$ is finite.
On the other hand, we have $c^2 (\delta_gx \otimes \delta_hy) =
Q^{x,y}_{g,h}Q^{y, x}_{h, g} \delta_hy  \otimes \delta_gx$, for all $g,
h \in G$, $x, y \in F$.
Thus the claim follows. \epf
\end{obs}

\begin{obs} Suppose that $R$ is a braided Hopf algebra. Then, by
\cite{T}, the antipode $\Ss$ commutes with $c$; that is, we have
$c(\Ss \otimes \id) = (\id \otimes \Ss)c$. This amounts to the
condition $Q^{x, y}_{g, h} = Q^{(g \fde x)^{-1}, y}_{(g \fiz
x)^{-1}, h}$, for all $x, y \in F$, $g, h \in G$, which
corresponds to the following relationship between  $\sigma$ and
$\tau$:
\begin{align*}& \sigma_{(h \fiz (g \fde x)^{-1})g}\left( x, (g \fiz
x)^{-1} \fde y \right) \,
\sigma_{h \fiz (g \fde x)^{-1}}(g \fde x, y)^{-1} \, \sigma_g(x,
(g \fiz x)^{-1} \fde y)^{-1} \\ & \tau_{x ((g \fiz x)^{-1} \fde
y)}(h \fiz (g \fde x)^{-1}, g) \, \tau_{(g \fiz x)^{-1} \fde y}
(h, g \fiz x)^{-1} \, \tau_x(h \fiz (g \fde x)^{-1}, g)^{-1} \\ &
= \sigma_{(h \fiz x)(g \fiz x)^{-1}}\left((g \fde x)^{-1}, g \fde
y \right) \, \tau_{(g \fde x)^{-1}(g \fde y)}(h \fiz x, (g \fiz
x)^{-1}) \\ & \sigma_{h \fiz x} (x^{-1}, y)^{-1} \, \tau_{g \fde
y}(g, g^{-1})^{-1} \, \sigma_{(g \fiz x)^{-1}} ((g \fde x)^{-1}, g
\fde y )^{-1} \, \tau_{(g \fde x)^{-1}} (h \fiz x, (g \fiz
y)^{-1}),
\end{align*}
for all $x, y \in F$, $g, h \in G$. \end{obs}

\bigbreak We now give a cohomological interpretation of
Proposition \ref{cond-q}.

Consider the following double complex:
\begin{equation}\label{double-comp}
{\bf C}^{..} : =
\begin{CD}
\vdots @.  \vdots \\ @AAA    @AAA \\ \Map_+(G^2 \times F,
\ku^{\times}) @>{\delta}>> \Map_+(G^2 \times F^2, \ku^{\times})
@>>> \dots
\\ @A{\delta'}AA    @A{\delta'}AA   \\ \Map_+(G \times F,
\ku^{\times})  @>{\delta}>> \Map_+(G \times F^2, \ku^{\times})
@>>> \dots,
\end{CD}
\end{equation}
where, for all $n, m \geq 1$, $\Map_+(G^n \times F^m,
\ku^{\times})$ is the abelian group of $\ku^{\times}$-valued
functions $f$ on $G^n \times F^m$ with the property that $f(g_n,
\dots, g_1; x_1, \dots, x_m) = 1$, if either one of $g_1, \dots,
g_n$ or $x_1, \dots, x_m$ is equal to $1$, and the maps $\delta$
and $\delta'$ are defined by
\begin{align*}
\delta f (g_q, \dots, g_1;& x_1, \dots, x_{p+1}) =
 f\left( g_q \fiz (g_{q-1} \dots g_1 \fde x_1), \dots, g_2 \fiz (g_1
\fde x_1), g_1 \fiz x_1; x_2, \dots, x_{p+1} \right) \\
& \times \prod_{i = 1}^p f(g_q \dots, g_1; x_1, \dots, x_ix_{i+1},
\dots, x_{p+1})^{(-1)^i}
 \times f(g_q, \dots, g_1; x_1, \dots, x_p)^{(-1)^{p+1}};
\end{align*}

\begin{align*}
\delta' f (g_{q+1}, \dots, g_1;& x_1, \dots, x_p)^{(-1)^p} =
 f\left( g_{q+1}, \dots, g_2; g_1 \fde x_1, (g_1 \fiz x_1) \fde x_2,
\dots, (g_1 \fiz x_1 \dots x_{p-1}) \fde x_p \right) \\
& \times \prod_{i = 1}^q f(g_{q+1}, \dots, g_{i+1}g_i, \dots, g_1; x_1,
\dots  x_p)^{(-1)^i}
 \times f(g_q, \dots, g_1; x_1, \dots, x_p)^{(-1)^{q+1}}.
\end{align*}
It is known  that the necessary and sufficient condition for $R$
to be a (usual) Hopf algebra is that the pair $(\sigma, \tau)$ be
a $1$-cocycle in the total complex $\Tot({\bf C}^{..})$, and
moreover that the assignment $(\sigma, \tau) \mapsto R$ defines an
isomorphism $H^1(\Tot ({\bf C}^{..})) \simeq \Opext (\ku^G, \ku
F)$. See \cite[Proposition 5.2]{Maext}.

Let $\delta^{\tot}$ denote the coboundary map in  the total
complex $\Tot({\bf C}^{..})$ and let
\begin{equation*}p_i:
\Tot({\bf C}^{..})^n = \oplus_{i = 1}^n\Map_+(G^{n-i} \times F^i,
\ku^{\times}) \to \Map_+(G^{n-i} \times F^i, \ku^{\times}),
\end{equation*}
 be the projection on the $i$-th. coordinate.
Conditions \eqref{norm-sigma}, \eqref{norm-tau},
\eqref{norm2-sigma} and \eqref{norm2-tau} say that $(\tau,
\sigma)$ belong to $\Map_+(G^{2} \times F, \ku^{\times}) \oplus
\Map_+(G \times F^2, \ku^{\times})$, and conditions
\eqref{cociclo-sigma}, \eqref{cociclo-tau} amount to
$p_1(\delta^{\tot}(\tau, \sigma)) = p_3(\delta^{\tot}(\tau,
\sigma)) = 1$.

\begin{cor} Let $c: R \otimes R \to R \otimes R$ be as in
\eqref{def-c}.
Then $R$ is a  prebraided Hopf algebra if and only if
\begin{equation}
Q_{g, h}^{x, y} = p_2(\delta^{\tot}(\tau, \sigma)) (g, h \fiz (g
\fde x)^{-1}; x, (g \fiz x)^{-1} \fde y) = \delta^{\tot}(\tau,
\sigma) (g, h \fiz (g \fde x)^{-1}; x, (g \fiz x)^{-1} \fde y),
\end{equation} for all $g, h \in G$,  $x, y \in F$. \end{cor}

\pf This is a reformulation of Proposition \ref{cond-q}. \epf

\begin{obs} It is natural to consider the following question:
given $Q \in \Map_+(G^2 \times F^2, \ku^{\times})$, find all pairs
of normalized cocycles $(\tau, \sigma)$ such that $R$ is a
prebraided Hopf algebra with braiding determined by $Q$ as in
\eqref{def-c}. This (possibly empty) space is a torsor over
$\Opext (\ku^G, \ku F)$. \end{obs}

\subsection{Braided compatible data for trivial actions.}

Along this subsection we shall assume that both actions $\fde$ and
$\fiz$ are trivial; that is, $\Sigma = F \times G$.
We discuss the compatibility conditions on the
cocycles $\sigma$ and $\tau$ in order that the corresponding
bicrossed product $R$ is a braided Hopf algebra  with non-trivial
braiding.

\bigbreak Let $\sigma: F \times F \to (\ku^{\times})^G$ and $\tau:
G \times G \to (\ku^{\times})^F$ be 2-cocycles satisfying the
normalization conditions \eqref{norm-sigma}, \eqref{norm-tau},
\eqref{norm2-sigma} and \eqref{norm2-tau}. We keep the notation
and conventions in \ref{braiding}.

By triviality of $\fde$ and $\fiz$, we may regard $\sigma$ and $\tau$ as
normalized maps $\sigma: G \to Z^2_+(F, \ku^{\times})$, $\tau: F
\to Z^2_+(G, \ku^{\times})$, and as such we may take their
differentials $\partial \sigma \in Z^2_+(G, Z^2_+(F,
\ku^{\times}))$, and $\partial \tau \in Z^2_+(F, Z^2_+(G,
\ku^{\times}))$.

\begin{lema}\label{act-triv} (i) Let $x, y \in F$, $g, h \in G$.
We have \begin{equation}Q^{x, y}_{g, h} = \partial(\sigma)(h,
g)(x, y)\, \, \partial(\tau)(x, y)(g, h). \end{equation}

(ii) Suppose that $\tau: F \to Z^2_+(G, \ku^{\times})$ is a group
homomorphism. Then $R$ is a braided Hopf algebra if and only if
$\partial \sigma \in \Hom (G/[G, G] \otimes G/[G, G], \Hom (F/[F,
F] \otimes F/[F, F], \ku^{\times}))$.

In this case, the braiding $c: R \otimes R \to R \otimes R$ is
trivial if and only if $\sigma$ is a group homomorphism.
\end{lema}

\pf Part (i) is straightforward. If $\tau$ is a homomorphism, then
$Q^{x, y}_{g, h} = \partial(\sigma)(h, g)(x, y)$. Thus, part (ii)
is a consequence of (i) and Lemma \ref{1-cociclos}. Clearly, $c$
is trivial if and only if $Q = 1$, if and only $\sigma$ is a group
homomorphism. \epf

\begin{obs} Let us fix $\tau: F \to Z^2_+(G, \ku^{\times})$ a group
homomorphism. If we start by taking $\sigma: G \to \Hom (F/[F, F]
\otimes F/[F, F], \ku^{\times})$, then $\sigma_g$ will be a
2-cocycle  on $F$, for all $g \in G$ (every bicharacter is), and the image of
$\partial \sigma$ will be contained in $\Hom (F/[F, F] \otimes F/[F, F],
\ku^{\times})$.
In order that the data be compatible, but $c$ be not trivial, we
need to have that $\partial \sigma$ is a bicharacter, but $\sigma$
is {\it not} a group homomorphism.  \end{obs}

\bigbreak Let now $\ku$ be the field $\mathbb C$ of complex
numbers. Let $p$ be an odd prime number and let $G = F = \F_p
\oplus \F_p$ be 2-dimensional vector spaces over the field $\F_p$
with $p$ elements. We use additive notation in both $G$ and $F$.
The elements of $F$ will be denoted with roman letters $(x, y),
(x', y'), \dots$, and the elements of $G$ will be denoted with
greek letters $(\alpha, \beta), (\alpha', \beta'), \dots$.

\bigbreak  We have $H^2(F, \mathbb C^{\times}) \simeq \Hom
(\Lambda^2F, \mathbb C^{\times})$ is of order $p$, and it consists
of the classes of the cocycles
\begin{equation*}\left( (x, y), (x', y') \right) \mapsto
\exp {\dfrac{2\pi i n}{p}(xy'-x'y)},
\end{equation*} where $n$ runs over the integers modulo $p$.

\begin{prop}\label{trivial-actions} Let $a$ and $b$ be integers modulo $p$ such that $ab \neq
0 \mod p$. Let $\tau \in \Hom (F, Z_+^2(G, \mathbb C^{\times}))$
be given by
\begin{equation}\tau_{(x, y)}\left( (\alpha, \beta), (\alpha', \beta')
\right)= \exp \left( \dfrac{2\pi i}{p} (x+y) (\alpha \beta' -
\alpha' \beta) \right), \qquad x, y, \alpha, \beta, \alpha',
\beta'  \in \F_p,
\end{equation}
and let $\sigma: G \to \Hom (\Lambda^2 F, \mathbb C^{\times})$ be
given by
\begin{equation}\sigma_{(\alpha, \beta)}\left( (x, y), (x', y')
\right) = \exp \left( \dfrac{2\pi i}{p} (a \alpha^2 + b \beta^2)
(xy'-x'y) \right), \qquad x, y, x', y', \alpha, \beta  \in \F_p.
\end{equation}
Then the associated bicrossed product $R$ is a braided Hopf
algebra, with non-trivial braiding $c$ given by \eqref{def-c}, and
\begin{equation}Q^{(x, y), (x', y')}_{(\alpha, \beta), (\alpha',
\beta')} = \exp \left( \dfrac{4\pi i}{p} (a \, \alpha \alpha' + b
\, \beta \beta') (xy'-x'y) \right),
\end{equation}
for all $x, y, x', y', \alpha, \beta,  \alpha', \beta'  \in \F_p$.

Moreover, $R$ is not commutative and not cocommutative. \end{prop}

\pf The proof follows from Lemma \ref{act-triv}. Note that, for
instance, $Q^{(1, 0), (0, 1)}_{(1, 0), (1, 0)} = \exp (\dfrac{4\pi
i}{p}a)$, which is not equal to $1$ since $p$ is odd and $a \neq 0
\mod p$. Finally, it is not difficult to see that $R$ is not
commutative and not cocommutative. \epf

\begin{obs} Observe that $Q = \partial \sigma$ is the symmetric
bilinear map associated with the quadratic map $\sigma$. Also,
$\tau$ is obtained as the composition of the epimorphism $F \to
\F_p$, $(x, y) \mapsto x+y$, with the natural isomorphism $(\F_p,
+) \simeq \Hom (\Lambda^2 F, \mathbb C^{\times})$. \end{obs}

\subsection{A categorical exact sequence.}
We shall assume in this subsection that $(R, c)$ is a braided Hopf
algebra. By Proposition  \ref{const-struct}, the braiding $c$ is
as described in Proposition \ref{const-struct}, where in addition,
the conditions in Theorem \ref{cond-sigma-tau} are satisfied.

\bigbreak Let $(V, c)$ be a braided vector space. Recall \cite{T}
that a subspace $W$ of $V$ is called  {\it categorical} if $c(V
\otimes W) \subseteq W \otimes V$ and $c(W \otimes V) \subseteq V
\otimes W$. In particular, if $W$ is a categorical subspace, then
it is a braided subspace with respect to $c\vert_{W \otimes W}: W
\otimes W \to W \otimes W$.

A quotient space $p: V \to U$ will be called {\it categorical} if
the kernel of $p$ is a categorical subspace of $V$. In this case,
$U$ is a quotient braided space with respect to the braiding $(p
\otimes p)(c): U \otimes U \to U \otimes U$.

\bigbreak The definition of extension of Hopf algebras may be
generalized to braided Hopf algebras as follows. We shall say that
the sequence of braided Hopf algebras and  braided Hopf algebra
maps $$ \begin{CD}1 @>>> S @> \iota>> R @>\pi>> T @>>> 1\end{CD}$$
is an \emph{extension of braided Hopf algebras} if $\iota$ is
injective, $\pi$ is surjective, $\ker \pi = RS^+$, $R^{\text{co }
\pi } = S$. We shall say that the extension is \emph{cleft} if
$\pi$ admits a section which is convolution invertible and
$T$-colinear.

\bigbreak Let $\iota: \ku^G \to R$ be the natural inclusion and
let $\pi: R \to \ku F$ be the natural projection.

\begin{prop}\label{categ-ext} The inclusion $\iota: \ku^G \to R$ and
the projection
$\pi: R \to \ku F$ are categorical.  Moreover, there is an exact
sequence of braided Hopf algebras
\begin{equation}\label{ex-sec}
\begin{CD}1 @>>> \ku^G @> \iota>> R @>\pi>> \ku F @>>>
1\end{CD}
\end{equation}
where the braiding in $\ku^G$ and $\ku F$ is the usual flip.
 \end{prop}

\pf The kernel of  $\pi$ is equal to the span of all elements
$\delta_gx$, where $g \in G \backslash \{ 1 \}$, and $x \in F$.
Note that since the braiding $c$  is diagonal in the basis
$\delta_gx$, $g \in G$, $x \in F$, it follows that $\iota$ and
$\pi$ are categorical.

By construction, $\iota$ is an algebra inclusion and $\pi$ is a
coalgebra surjection. Also, it is not difficult to see that
condition \eqref{norm2-sigma} is equivalent to $\iota: \ku^G \to
R$ being a coalgebra map, while condition \eqref{norm2-tau} is
equivalent to $\pi: R \to \ku F$ being an algebra map. The rest of
the proposition follows easily. \epf

\begin{obs} $R^*$ is also a braided Hopf algebra, which can be
constructed from the matched pair arising from the exact
factorization $\Sigma =GF$, and it fits into an exact sequence of
braided Hopf algebras $\begin{CD}1 @>>> \ku^F @> \pi^*>> R^*
@>\iota^*>> \ku G @>>> 1\end{CD}$.
\end{obs}

\begin{definition}We shall say that the extension of braided Hopf algebras
\eqref{ex-sec}  is  \emph{realizable} over $H$, whenever
the braided Hopf algebras $\ku^G$, $R$ and $\ku F$, as well
as the maps $\iota$ and $\pi$, are in the category
${}^H_H\mathcal{YD}$.
\end{definition}

 It follows from Proposition \ref{categ-ext} and the results
in \cite[Section 6]{T} that there exists a Hopf algebra $H$ such
that \eqref{ex-sec} is realizable over $H$. Indeed, there is a Hopf algebra $H$ such that
$\iota$ is categorical by \cite[Proposition 6.6]{T}; but then $\pi$ is also categorical.

\bigbreak But it is \emph{not} true that if $R$ is realizable over
any Hopf algebra $K$ then $\ku^G$ and $\ku F$ also are. For
instance, assume that $\ku$ is algebraically closed of
characteristic zero and let $R$ be the group algebra of a finite
group $L$, $N$ a normal abelian subgroup, $G$ the group of
characters of $N$ and $F = L/N$. Let $\theta$ be a non-trivial
automorphism of $L$, say of finite order, and let $C$ be the
subgroup of $\Aut L$ generated by $\theta$. Then $R$ is a
Yetter-Drinfeld module over $\ku C$ with trivial coaction, but
$\ku^G$ is not a Yetter-Drinfeld submodule of $R$ unless $N$ is
$\theta$-stable. For a concrete example, let $N$ be a finite
abelian group, $L = N \times N$ and $\theta$ the transposition.

\bigbreak
Let $S$ be any braided Hopf algebra.
Consider the left {\it braided adjoint action} of $S$ on itself,
given by
\begin{equation}
\ad_c(a)(b) = m \left( \id \otimes m (\id \otimes \Ss) c \right)
(\Delta \otimes \id) (a \otimes b), \qquad a, b \in S.
\end{equation}

Let $H$ be a Hopf algebra such that $S$ is realizable over $H$.
Then the left braided adjoint action $\ad_c$ of $S$ coincides with
the restriction to $S$ of the left adjoint action of the
corresponding Radford biproduct $S \# H$: $\ad_c(a) (b) = \ad_{S
\# H}(a) (b)$, for all $a, b \in S$.

We come back to our situation.  The next lemma shows that $\ku^G$
is 'braided normal' in $R$.

\begin{lema} We have $\ad_c(\delta_gx)
(\delta_h) = \delta_{g, 1} \delta_{(h^{-1} \fiz x^{-1})^{-1}}$,
 for all $g, h \in G$, $x \in F$.

In particular, the categorical braided  Hopf subalgebra $\ku^G
\subseteq R$ is stable under the left braided adjoint action.
\end{lema}

\pf Straightforward. \epf

\subsection{Equivalences}
Let $R$ and $R'$ be braided Hopf algebras. A linear map $\Theta: R
\to R'$ is called a \emph{morphism of braided Hopf algebras} if it
preserves the multiplication, comultiplication, unit and counit
maps. Since the antipode is the convolution inverse of the
identity, it follows that any morphism $\Theta$ of braided Hopf
algebras preserves also the antipode. Hence, by \cite{Sb},
$\Theta$ commutes with the braiding; that is, $(\Theta \otimes
\Theta) c_R = c_{R'}(\Theta \otimes \Theta)$.

\bigbreak Let $\begin{CD}1 @>>> S @> \iota>> R @>\pi>> T @>>>
1\end{CD}$ and $\begin{CD}1 @>>> S @> \iota'>> R' @>\pi'>> T @>>>
1\end{CD}$ be two extensions of braided Hopf algebras. An
isomorphism $\Theta : R \to R'$ of braided Hopf algebras is an
\emph{isomorphism of extensions} if the following diagram commutes
\begin{equation*}\begin{CD}1 @>>> S @> \iota>> R @>\pi>> T @>>>
1 \\ @. @V{\id}VV  @V{\Theta}VV   @V{\id}VV @. \\ 1 @>>> S @>
\iota'>> R' @>\pi'>> T @>>> 1. \end{CD}
\end{equation*}

\begin{prop}\label{prop-equiv} Let $R = \ku^G {}^{\tau}\#_{\sigma}\ku F$ and
$R'  =  \ku^G {}^{\tau'}\#_{\sigma'}\ku F$ be braided Hopf
algebras and consider the corresponding extensions as in
\eqref{ex-sec}. Let $\nu \in \Map_+(G \times F, \ku^{\times})$ and
define $\Theta: R \to R'$ in the form $\Theta (\delta_gx) = \nu(g,
x) \delta_gx$, for all $g \in G$, $x \in F$. Then $\Theta$ is an
isomorphism of extensions if and only if $(\tau, \sigma) = (\tau',
\sigma') \, \delta^{\tot}\nu$ in the complex \eqref{complejo}.

Furthermore, any isomorphism of extensions $\Theta: R \to R'$
arises in this way for a unique $\nu$.
\end{prop}

\pf It is easy to see that $\Theta$ is an algebra map if and only
if $\sigma_g(x, y) \nu (g, xy) = \sigma'_g(x, y) \nu(g, x) \nu(g
\fiz x, y)$, and $\nu(g, 1) = 1$, for all $g \in G$, $x, y \in F$.
Also, $\Theta$ is a coalgebra map if and only if $\tau'_x(g, h)
\nu (gh, x) = \tau_x(g, h) \nu(g, h \fde x) \nu(h, x )$, and
$\nu(1, x) = 1$, for all $g, h \in G$, $x \in F$. This proves the
first claim.

Let now $\Theta: R \to R'$ be an isomorphism of extensions of
braided Hopf algebras. Since $\pi = \pi' \Theta$, it can be seen
that $\Theta(x)x^{-1} \in \ku^G$, for all $x \in F$. Define
$\nu(g, x)$ by $\Theta(x)x^{-1} = \sum_{g \in G}\nu(g, x)
\delta_g$. Then $\Theta (\delta_g x) = \delta_g \Theta (x) =
\nu(g, x) \delta_g x$, and the conclusion follows from the first
claim.
 \epf

\begin{cor} The group of automorphisms of the extension
\eqref{ex-sec}is isomorphic to  $Z^0(\Tot (\mathbf C^{..}))$. Any
such automorphism is categorical. \qed \end{cor}

\begin{obs}\label{sigma-coborde} Suppose that the 2-cocycle $\sigma: F \times F \to (\ku^{\times})^G$
is a coboundary. Then $R$ is isomorphic to a bicrossed product
$\ku^G \, {}^{\tau'}\# \ku F$. In particular, this always happens
if all the Sylow subgroups of  $F$ are cyclic; see \cite[Lemma
1.2.5]{pqq}.
\end{obs}

\subsection{Commutativity} We shall say that a braided Hopf algebra $R$ is
 \emph{braided commutative} if $m = m c: R \otimes R \to R$; respectively, $R$ is
called \emph{braided cocommutative} if $\Delta = c \Delta : R \to
R \otimes R$.

\bigbreak Let $R = \ku^G \, {}^{\tau}\#_{\sigma} \ku F$. The
verification of the following claims is straightforward:

\bigbreak (a) $R$ is braided commutative if and only if $F$ is
abelian, $\fiz$ is trivial and
\begin{equation}\label{br-comm}Q^{x, y}_{g, g} = \sigma_g (x, y)
\sigma_g (y, x)^{-1}, \qquad \forall x, y \in F, \, g \in G.
\end{equation}

\bigbreak (b) $R$ is braided cocommutative if and only if $G$ is
abelian, $\fde$ is trivial and
\begin{equation}\label{br-cocomm}Q^{x, x}_{g, h} = \tau_x (h, g)
\tau_x (g, h)^{-1}, \qquad \forall x \in F, \, g, h \in G.
\end{equation}

\section{Diagonal realizations over finite groups}\label{diag-real}
We shall consider in this section $2$-cocycles $\sigma: F \times F
\to (\ku^{\times})^G$ and $\tau: G \times G \to (\ku^{\times})^F$
satisfying the  normalization conditions \eqref{norm-sigma},
\eqref{norm-tau}, \eqref{norm2-sigma} and \eqref{norm2-tau}. We
discuss a particular but important class of realizations.

\bigbreak We fix a finite group $C$ and we let $H = \ku C$.
 We fix functions $z: G \times F\to Z(C)$ and $\chi: G \times F\to
 \wC$,
and we define a structure of left Yetter-Drinfeld module on $R$ by
imposing $\delta_g x \in R^{\chi(g,x)}_{z(g,x)}$, $g\in G$, $x\in
F$. That is, the action and the coaction of $H$ on $R$ are given,
respectively, by
\begin{equation} u. \delta_gx = \la \chi(g, x), u\ra \, \delta_gx,
\qquad u\in C; \qquad \rho (\delta_gx) = z(g, x) \otimes
\delta_gx.\end{equation} In particular, the braiding $c: R \otimes
R \to R \otimes R$  is given in this case by
\begin{equation}\label{tak-comp}c (\delta_gx \otimes \delta_hy) = \la \chi(h, y),
z(g, x)\ra \, \delta_hy \otimes \delta_gx.\end{equation}

\begin{lema}\label{uno-cociclos} (i). The multiplication of $R$ given
by
\eqref{producto-def}
is a morphism of $H$-modules if and only if
\begin{equation}
\label{uno-cociclo1} \chi(g, xy) =  \chi(g, x) \chi(g \fiz x, y),
\qquad g\in G, x, y\in F.
\end{equation}
(ii).  The comultiplication of $R$ given by \eqref{coproducto-def}
is a morphism of $H$-modules if and only if
\begin{equation}
\label{uno-cociclo2} \chi(gh, x) =  \chi(g, h\fde x) \chi(h, x),
\qquad g, h\in G, x\in F.
\end{equation}
(iii).  The multiplication of $R$
is a morphism of $H$-comodules if and only if
\begin{equation}
\label{uno-cociclo3} z(g, xy) =  z(g, x) z(g \fiz x, y), \qquad
g\in G, x, y\in F.
\end{equation}
(iv).  The comultiplication of $R$
is a morphism of $H$-comodules if and only if
\begin{equation}
\label{uno-cociclo4} z(gh, x) =  z(g, h\fde x) z(h, x), \qquad g,
h\in G, x\in F.
\end{equation}
\end{lema}
\pf Straightforward. \epf

\begin{obs}\label{norm-c-chi} The following normalization properties
follow from \eqref{uno-cociclo1}--\eqref{uno-cociclo4}:
\begin{equation}\label{n-cchi} z(1, x) = 1, \qquad z(g, 1) = 1,  \qquad
\chi(1, x) = 1,
\qquad \chi(g, 1) = 1, \qquad x \in F,\, g \in G. \end{equation}
It is not difficult to see that these amount to the unit and
counit maps being  morphisms of Yetter-Drinfeld modules.
\end{obs}

\begin{exa} The action $\fde$ induces a right action of $G$ on $\Map_+(F,
Z(C))$ in the form $(\phi \leftharpoonup g) (x) = f(g \fde x)$.
Let $\psi \in \Map_+(F, Z(C))$ and consider the function
$z_{\psi}: G \times F \to Z(C)$ given by $z_{\psi}(g, x) =
\partial_{\psi} (g) (x) = \psi(g \fde x) \psi^{-1}(x)$. Then $z_{\psi}$
satisfies the cocycle condition \eqref{uno-cociclo4} by
construction: indeed, $z_{\psi}$ is the $1$-coboundary of $\psi$
'in the first variable'.

\begin{lema} A sufficient condition for
$z_{\psi}$ to satisfy condition \eqref{uno-cociclo3}  is that
$\psi: F \to Z(C)$ be a group homomorphism.
\end{lema}

\pf We compute, for all $x, y \in F$, $g \in G$,
\begin{align*}
z_{\psi}(g, xy) & = \psi(g \fde (xy)) \psi^{-1}(xy),\\ z_{\psi}(g,
x) z_{\psi}(g \fiz x, y) & = \psi(g \fde x) \psi((g \fiz x) \fde
y)) \psi^{-1}(x) \psi^{-1}(y);
\end{align*}
In view of \eqref{comp1}, both expressions are equal whenever
$\psi$ is a group homomorphism. \epf
\end{exa}

\bigbreak
We give now an alternative approach to Theorem \ref{cond-sigma-tau}.
The following theorem is a consequence of Proposition
\ref{cond-q}.

\begin{theorem}\label{braid-c-chi} Suppose that $z: G \times F \to
Z(C)$ and
$\chi: G \times F \to \widehat C$ satisfy conditions
\eqref{uno-cociclo1}--\eqref{uno-cociclo4} in Lemma
\ref{uno-cociclos}. Then $R$ is a braided Hopf algebra over $\ku
C$ if and only if
\begin{equation}\label{comp-cchi-2coc}
\sigma_{ts}(x, y) \tau_{xy}(t, s)  = \la \chi \left(t \fiz (s \fde
x), (s \fiz x) \fde y\right), z(s, x) \ra \,
 \tau_x(t, s) \, \tau_y(t \fiz (s \fde x), s \fiz x) \,
\sigma_{t}(s \fde x, (s \fiz x) \fde y) \, \sigma_{s}(x, y),
\end{equation}
for all $s, t \in G$, $x, y \in F$. If this holds, we shall say
that $(z, \chi)$ is a \emph{diagonal realization} of $R$ over $\ku
C$. \qed
\end{theorem}

\begin{obs}\label{naif} Consider the conditions
\begin{align}\label{comp-cchi-sigma}
\sigma_{ts}(x, y)   & = \la \chi \left(t \fiz (s \fde x), (s \fiz
x) \fde y\right), z(s, x) \ra \, \sigma_{t}(s \fde x, (s \fiz x)
\fde y) \, \sigma_{s}(x, y),  \\ \label{comp-cchi-tau}
\tau_{xy}(t, s)  & =
 \tau_x(t, s) \, \tau_y(t \fiz (s \fde x), s \fiz x),
\end{align}
for all $s, t \in G$, $x, y \in F$. It is clear that any two among
\eqref{comp-cchi-2coc}, \eqref{comp-cchi-sigma} and
\eqref{comp-cchi-tau} imply the third. This observation will be
used later in order to systematically produce examples of diagonal
realizations.

A similar observation applies if one considers instead the
conditions
\begin{align}\label{comp-cchi-sigma2}
\sigma_{ts}(x, y)   & =  \sigma_{t}(s \fde x, (s \fiz x) \fde y)
\, \sigma_{s}(x, y),  \\ \label{comp-cchi-tau2} \tau_{xy}(t, s)  &
= \la \chi \left(t \fiz (s \fde x), (s \fiz x) \fde y\right), z(s,
x) \ra \, \tau_x(t, s) \, \tau_y(t \fiz (s \fde x), s \fiz x).
\end{align}

\bigbreak Suppose that $R = \ku^G {}^{\tau}\#_{\sigma} \ku F$ is a
braided Hopf algebra and the 2-cocycle $\sigma: F \times F \to
(\ku^{\times})^G$ is a coboundary. In view of Remark
\ref{sigma-coborde}, $R$ is isomorphic to a  bicrossed product
$\ku^G {}^{\tau'}\# \ku F$. Since $\sigma' = 1$ satisfies
\eqref{comp-cchi-sigma2}, then $\tau'$ must satisfy
\eqref{comp-cchi-tau2}.

This simplifies the search of braided compatible  data in many
cases, for instance, in the case where all Sylow subgroups of $F$
are cyclic.
\end{obs}

\bigbreak Suppose that $R$ admits a diagonal realization over $\ku
C$, and consider the Radford biproduct $R \# \ku C$.

\begin{prop} (i) The extension
 $\begin{CD} 1 @>>> \ku^G @>{\iota}>> R @>{\pi}>> \ku F @>>>
1\end{CD}$ is realizable over $\ku C$;

(ii) there are exact sequences of Hopf algebras
\begin{align}\label{ex-ses1}
& 1 \to \ku^G \to R \# \ku C \to \ku F \otimes \ku C \to 1, \\
\label{ex-ses2} & 1 \to \ku^G  \otimes \ku C \to R \# \ku C \to
\ku F \to 1,
\end{align}
where all maps are canonical. \end{prop}

\pf (i). Consider the trivial action and coaction  of $\ku C$ on
$\ku^G$ and $\ku F$, making them Yetter-Drinfeld modules. It
follows from the normalization conditions \eqref{n-cchi}, that the
canonical maps $\iota$ and $\pi$ are morphisms of Yetter-Drinfeld
modules. This proves (i).

Note that for the corresponding biproducts we have  $\ku^G \# \ku
C = \ku^G  \otimes \ku C$, and $\ku F \# \ku C = \ku F  \otimes
\ku C$

(ii). Conditions \eqref{n-cchi} also imply that the action and
coaction of $\ku C$ on  elements $\delta_g$,  $g \in G$, and also
on elements $x \in F$, are both trivial. Using this plus part (i),
one sees that the maps in \eqref{ex-ses1} and \eqref{ex-ses2} are
Hopf algebra maps. The exactness follows easily. \epf

\bigbreak We now give an interpretation of the conditions in Lemma
\ref{uno-cociclos} in the terms of the  cohomology of a complex
closely related to that considered in \eqref{double-comp}.  Let
$M$ be an abelian group. Let $n, m \geq 1$, and let $\Map_+(G^n
\times F^m, M)$ be the abelian group of $M$-valued functions $f$
on $G^n \times F^m$ with the property that $f(g_n, \dots, g_1;
x_1, \dots, x_m) = 1$, if either one of $g_1, \dots, g_n$ or $x_1,
\dots, x_m$ is equal to $1$. Consider the double complex
\begin{equation}\label{complejo}
{\bf C}^{..} (M): =
\begin{CD}
\vdots @.  \vdots \\
@AAA    @AAA \\
\Map_+(G^2 \times F, M) @>{\delta}>> \Map_+(G^2 \times F^2, M) @>>>
\dots \\
@A{\delta'}AA    @A{\delta'}AA   \\
\Map_+(G \times F, M)  @>{\delta}>> \Map_+(G \times F^2, M) @>>> \dots,
\end{CD}
\end{equation}
where the maps $\delta$ and $\delta'$ are defined as for the
complex \eqref{double-comp}.

\begin{lema} (i) Conditions \eqref{uno-cociclo1} and
\eqref{uno-cociclo2} are equivalent to $\chi \in Z^0(\Tot ({\bf
C}^{..} (\widehat C)))$.

(ii) Conditions \eqref{uno-cociclo3} and \eqref{uno-cociclo4} are
equivalent to $c \in Z^0(\Tot ({\bf C}^{..} (Z(C))))$. \end{lema}

\pf Straightforward. \epf

We close this section by showing that diagonal realizations over
abelian groups always exist. Suppose that $R$ is a braided Hopf
algebra; so that the conditions in Theorem \ref{cond-sigma-tau}
are satisfied. Let $Q^{x,y}_{g,h}$ be given by formula
\eqref{def-q}; thus $Q$ satisfies the conditions in Lemma
\ref{1-cociclos}. Let $\Gamma$ be either $\Z$, or $\Z/N$ provided
that the order of $Q^{x,y}_{g,h}$ divides $N$ for all $x,y \in F$,
$g,h\in G$. Consider the group $\Gamma^{G\times F}$, with the
canonical elements $e(g,x)$. We then define
\begin{align*}
C &:= \Gamma^{G\times F} / \la \{e(g, xy) - e(g,x) - e(g\fiz x,
y); \, e(gh, x) - e(g, h \fde x) - e(h, x): g, h\in G, x,y \in
F\}\ra,
\\ z(g,x) &:= \text{ the class of } e(g,x) \text{ in } C.
\end{align*}

\begin{lema}\label{close} There are characters $\chi(g,x)$ of $C$ defined by
$$ \chi(g,x) \left(z(h,y)\right)
= Q^{x,y}_{g,h}.$$ Furthermore,
$(z, \chi)$ is a diagonal realization of $R$ over $\ku
C$. \end{lema}

\pf The characters are well-defined by \eqref{comm1-m} and \eqref{comm1-delta};
this is a diagonal realization by the definition of $z$, \eqref{comm2-m} and \eqref{comm2-delta}.
\epf

\section{Examples of diagonal realizations}\label{eg-diag-real}

We discuss in the next subsections, under additional assumptions
on the matched pair $(G, F)$, some reductions in order to
determine maps $c : G \times F \to Z(C)$ and $\chi : G \times F
\to \widehat C$ satisfying the conditions in Lemma
\ref{uno-cociclos}. We shall assume in this section that $\ku$ is
algebraically closed of characteristic zero.

\subsection{Semidirect products.}

Consider the case where the action $\fde$ is  trivial; so that
$\fiz : G \times F \to G$ is an action by group automorphisms and
the group $\Sigma = FG$ is isomorphic to the associated semidirect
product $F \ltimes G$.

The action $\fiz$ induces by transposition  left actions of $F$ on
$\Hom(G, Z(C))$ and on $\Hom(G, \widehat C)$; for instance, we
have $(x \rightharpoonup \phi) (g) = \phi(g \fiz x)$,  for all $x
\in F$, $g \in G$, $\phi \in \Hom(G, Z(C))$.

\begin{lema}\label{z1-g} (i) The set of maps $\chi: G \times F \to
\widehat C$
satisfying
\eqref{uno-cociclo1} and \eqref{uno-cociclo2} is in bijective
correspondence with $Z^1(F, \Hom(G, \widehat C))$.

(ii) The set of maps $z: G \times F \to Z(C)$ satisfying
\eqref{uno-cociclo3} and \eqref{uno-cociclo4} is in bijective
correspondence with $Z^1(F, \Hom(G, Z(C)))$.
\end{lema}

\pf We prove (i), the proof of (ii) being similar. The
correspondence is given by $\chi \mapsto \tilde \chi: F \to
\Hom(G, \widehat C)$, $\tilde \chi (x) (g) = \chi (g, x)$.
Condition \eqref{uno-cociclo2} amounts to $\tilde \chi (x) \in
\Hom(G, \widehat C)$, for all $x \in F$, and condition
\eqref{uno-cociclo1} says exactly that $\tilde \chi$ is a
$1$-cocycle. \epf

\begin{cor}\label{c=1} Suppose that $|G|$ and $|Z(C)|$ are relatively
prime. If
$z: G \times F \to Z(C)$ satisfies \eqref{uno-cociclo3} and
\eqref{uno-cociclo4}, then $z(g, x) = 1$, for all $x \in F$, $g
\in G$. \end{cor}

\pf In this case we have  $\Hom(G, Z(C)) = 1$. Therefore the claim
follows from Lemma \ref{z1-g}. \epf

\begin{exa} Let $p$ be a prime number and suppose that
$\dim R = p^3$. In other words, $\Sigma$ has order $p^3$.
Up to passing to the dual, we may assume that $|G| = p^2$ and
thus that $G$ is normal in $\Sigma$.
 Assume also that $p$
does not divide the order of $C$. Then Corollary \ref{c=1} implies
that $c = 1$. Thus $Q^{x, y}_{g, h} = 1$, for all $x, y \in F$,
$g, h \in G$. Hence every $R$ arising from this setup is trivial,
{\it i.e.}, is a usual Hopf algebra.
\end{exa}

\bigbreak
We now fix $\widetilde z \in Z^1(F,  \Hom(G, Z(C)))$
and $\widetilde \chi \in Z^1(F, \Hom(G, \wC))$; they are the 1-cocycles
corresponding to maps $z: G \times F \to Z(C)$ and $\chi: G \times
F \to \widehat C$, respectively, as in Lemma \ref{z1-g}.
We look for $\sigma$, $\tau$ satisfying the conditions in Theorem \ref{braid-c-chi}.

Consider the action of $F$ on $Z^2(G, \ku^{\times})$  given by
$(x.f) (g, h) = f(g \fiz x, h \fiz x)$, $g, h \in G$, $x \in F$;
this action is well-defined because $F$ acts by group
automorphisms on $G$.

The map $\tau: G \times G \to (\ku^{\times})^F$ can be regarded as
a map $F \to \Map (G \times G, \ku^{\times})$; we shall write
$\widetilde \tau$ to indicate this latter map. Note that $\tau$ is
a 2-cocycle if and only if the image of $\widetilde \tau$ is
contained in $Z^2(G, \ku^{\times})$.

\begin{prop}\label{gen-som} Let $\widetilde z \in Z^1(F, \Hom(G,
Z(C)))$
and $\widetilde \chi \in Z^1(F, \Hom(G, \wC))$.  Let also $\sigma
\in Z^2(F, (\ku^{\times})^G)$ be a normalized 2-cocycle such that
the following compatibility condition holds:
\begin{equation}\label{cond-som}
\sigma_{ts}(x, y)  = \la (x \rightharpoonup \widetilde{\chi})(y)
(s),  \widetilde z(x) (t) \ra \, \sigma_{s}(x, y) \, \sigma_{t}(x,
y).
\end{equation}
Let $\tau: G \times G \to (\ku^{\times})^F$ be a normalized
2-cocycle. Assume that the normalization conditions
\eqref{norm2-sigma} and \eqref{norm2-tau} hold.

Then $R$ is a braided Hopf algebra over $\ku C$ if and only if
$\widetilde \tau: F \to Z^2(G, \ku^{\times})$ is a 1-cocycle.
\end{prop}

\pf This is a special instance of Remark \ref{naif}. In view of
Theorem \ref{braid-c-chi}, and using \eqref{cond-som}, $R$ is a
braided Hopf algebra if and only if
\begin{equation*}\tau_{xy}(t, s) = \tau_x(t, s) \,
\tau_y(t \fiz x, s \fiz x) = \tau_x(t, s) \, (x . \tau_y)(t, s),
\end{equation*}
for all $s, t \in G$, $x, y \in F$, which is exactly the 1-cocycle
condition on $\widetilde \tau$ . This proves the proposition. \epf

\begin{exa}\label{nicosomm}
If $\tau$ is the trivial 2-cocycle, then $\widetilde \tau$ is the
trivial 1-cocycle, and we get a braided Hopf algebra $R$; note
that $R$ is the tensor product $\ku^G \otimes \ku F$ as a
coalgebra.  This braided Hopf algebra structure of $R$ is due to
Andruskiewitsch and Sommerh\" auser. Its construction appears in
\cite[3.2]{S}.
\end{exa}

The proposition above can be used to construct examples of
non-trivial braided  Hopf algebras, using the data in \cite[Ch.
3]{S}. These examples are not commutative and also not
cocommutative. A particular case of this construction is done in
the next subsection.

\subsection{An example from finite fields.}

We first recall the construction in \cite[3.3]{S}. Let $K$ be a
finite ring, let $F$ be a finite group, let $G$
denotes the additive group of $K$ written additively
and let there be given the following data:
\begin{flalign}
\label{action} & \text{a group homomorphism } \, \nu: F \to
K^{\times}&
\end{flalign}
(we endow $G$ with the $F$-action
defined by $g \fiz x : = g \nu (x)$, $x \in F$, $g \in G$),
\begin{flalign}\label{1-coc}
& \text{two 1-cocycles } \, \alpha, \beta \in Z^1(F,
G), &
\\ \label{2-coc} & \text{a normalized 2-cocycle } \, \phi \in Z^2(F,
G),
&
\\ \label{trazas} & \text{two characters } \, \eta, \lambda: G \to
\ku^{\times}, \, \text{such that} \, \la \lambda, ghs \ra = \la
\lambda, hgs \ra, \, g, h, s \in K. & \end{flalign}

$(F, G)$ is a matched pair with respect to the action
$\fiz: G \times F \to G$ and the trivial action $\fde : G \times F
\to F$. With respect to the above data, we define $z: G \times F
\to K$, $\chi: G \times F \to K$ and $\sigma: F \times F \to
(\ku^{\times})^G$ in the form
\begin{flalign}
\label{ex-c} & z(g, x) = g \beta(x), &
\\ \label{ex-chi} & \la \chi(g, x), h \ra = \la \lambda, h g
\alpha(x)\ra^2, &
\\ \label{ex-sigma} & \sigma_g(x, y) = \la \eta, g \phi(x, y) \ra \,
\la \lambda, g^2 \nu(x) \beta(x) \alpha(y) \ra, & \end{flalign}
for all $x, y \in F$, $g, h \in G$. The result in \cite[3.3]{S},
combined with Lemma \ref{gen-som}, implies that for any normalized
1-cocycle $\widetilde \tau \in Z^1(F, Z^2(G, \ku^{\times}))$, the
associated bicrossed product $R$ is a braided Hopf algebra over
$\ku C$, where $C$ is the additive group of $K$.

\bigbreak We construct now explicit examples of 1-cocycles
$\widetilde \tau$ in this situation. For notational simplicity, we
assume that $\ku$ is the field $\mathbb C$ of complex numbers. Let
$p$ and $q$ be prime numbers such that $p = 1 \mod q$,  let $K =
\F_{p^2}$ be the field with $p^2$ elements, and let $F = \Z_q$ be
the cyclic group of order $q$.

The assumption on $q$ implies that there exists a unit modulo $p$,
$\nu \in \F_p^{\times} \subseteq K^{\times}$, of order $q$; by
abuse of notation, we let $\nu: F \to K^{\times}$ be the group
homomorphism given by $\nu (x) = \nu^x$.

As before, let $G$ denote the additive group of $K$, and consider
the right action of $F$ on $G$ given by $g \fiz x = \nu^x g$.

\begin{lema}\label{q-no} Let $0 \neq \alpha: F \to G$ be a map.
Then $\alpha$ is a normalized 1-cocycle if and only of
there exists $r \in K^{\times}$ such that $\alpha$ has the form
\begin{equation}\alpha(x) = r[x]_{\nu}, \qquad x \in F - 0,\end{equation}
$\alpha(0) = 0$. \end{lema}

Here, $[x]_{\nu}$ denotes the $\nu$-number: $[x]_{\nu}: = 1 + v +
\dots + v^{x-1} \in \F_p$.

\pf The 1-cocycle condition on $\alpha$ says that $\alpha (x + y)
= \alpha(x) + \nu^x \alpha(y)$, for all $x, y \in F$. Using that
$F$ is cyclic, generated by $1$, one can show by induction that
$\alpha (x) = \alpha(1) [x]_\nu$, for all $x \in F$. Putting $r =
\alpha(1)$, which is non-zero by assumption, the lemma follows.
\epf

Let $a \in K$ such that $K = \F_p(a)$; so that every element $g
\in K$ writes uniquely in the form $g = j + la$, with $j, l \in
\F_p$. This determines an isomorphism between the additive group
of $K$ and $\F_p \oplus \F_p$. We shall denote $\det_a: K \times K
\to \F_p$ the function defined by $\det_a(g, h) = jl' - lj'$, for
$g = j + la$, $h =  j' + l'a \in K$.

\begin{lema}\label{def-tau} Let $x \in F$, $g, h \in G$. The formula
\begin{equation}\tau_x(g, h) =
\exp \left( \dfrac{2\pi i}{p} \, [x]_{\nu^2} \, {\det{}_a}(g, h)
\right),
\end{equation} defines a 1-cocycle $\widetilde \tau \in Z^1(F, Z^2(G, \mathbb
C^{\times}))$.
\end{lema}

\pf Every cohomology class in $H^2(G, \mathbb C^{\times}) \simeq
\F_p$ can be represented by one of the 2-cocycles $\kappa_n: (g,
h) \mapsto \exp {\dfrac{2\pi i}{p} n \, \det{}_a(g, h)}$, where $n$
runs over the integers modulo $p$, giving a group isomorphism
$\F_p \simeq H^2(G, \mathbb C^{\times})$. We have
\begin{align*}(x . \kappa_n) (g, h) & = \exp {\dfrac{2\pi i}{p} n \,
\det{}_a(g \fiz x,
h \fiz x)} = \exp {\dfrac{2\pi i}{p} n \, \det{}_a(g \nu^x, h
\nu^x)}
\\ & = \exp {\dfrac{2\pi i}{p} n \, \nu^{2x} \, \det{}_a(g , h)},
\end{align*}
the last equality because we have chosen $\nu \in \F_p$. Thus the
action of $x \in F$ on $\F_p \simeq H^2(G, \mathbb C^{\times})$ is
given by multiplication by $\nu^{2x}$ .

The argument in the proof of Lemma \ref{q-no} shows that every
1-cocycle $\widetilde \tau: F \to H^2(G, \mathbb C^{\times})$ is
of the form $\tau(x) = \kappa_{r [x]_{\nu^2}}$, for some $r \in
\F_p^{\times}$. This implies the lemma. \epf

\begin{prop}\label{eg-p4q} Let $C$ be the additive group of $\F_{p^2}$. Let $\nu \in
\F_p^{\times}$, and consider the matched pair $(F, G)$ as above.
Let also $z: G \times F \to C$, $\chi: G \times F \to C$, $\sigma:
F \times F \to (\ku^{\times})^G$ and $\tau: G \times G \to
(\ku^{\times})^F$, be defined by
\begin{flalign}
\label{exa-c} & z(g, x) = g [x]_{\nu}, &
\\ \label{exa-chi} & \la \chi(g, x), h \ra = \exp \left( \dfrac{4\pi
i}{p} \,
\tr (h g) \, [x]_{\nu} \right), &
\\ \label{exa-sigma} & \sigma_g(x, y) =
\exp \left( \dfrac{2\pi i}{p} \, \tr (g^2) \, \nu^x [x]_{\nu}
[y]_{\nu} \right), &
\\ & \tau_x(g, h) = \exp \left( \dfrac{2 \pi i}{p} \, [x]_{\nu^2}
\, \det{}_a (g, h) \right), &
\end{flalign}
for all $x, y \in F$, $g, h \in G$, where $\tr: \F_{p^2} \to \F_p$
is the trace map.

Then the associated bicrossed product $R$ is a braided Hopf
algebra over $\mathbb C C$. The braiding on $R$ is given by
\eqref{def-c}, where
\begin{equation}Q^{x,  y}_{g, h} = \exp \left( \dfrac{4 \pi i}{p} \,
\tr (h g) \, [x]_{\nu} \,
[y]_{\nu} \right), \qquad g, h \in \, x, y \in F.
\end{equation}
\end{prop}

Note that $R$ is non-trivial and also not commutative and not
cocommutative. The dimension of $R$ is $p^2q$, and the dimension
of the biproduct $R \# \mathbb C C$ is $p^4q$.

\pf The proof follows from \cite[3.3]{S} and Lemma \ref{gen-som},
using Lemma \ref{def-tau}. Take $\phi = 1$, $\eta = 1$; $\nu: F
\to G$, group homomorphism given by $\nu (x) = \nu^x$; 1-cocycles
$\alpha = \beta: F \to G$, with $\alpha (x) = [x]_{\nu}$; and let
$\lambda : G \to \mathbb C^{\times}$ be the group homomorphism
defined by $\la \lambda, g \ra = \exp (\dfrac{2 \pi i}{p} \, \tr
(g))$, for all $x \in F$, $g \in G$. \epf

\subsection{Direct products} Suppose that {\it both} actions
$\fde$ and $\fiz$ are  trivial, {\it i.e.}, that  $\Sigma \simeq F
\times G$. In this case, the maps  $z: G \times F \to Z(C)$
satisfying \eqref{uno-cociclo3} and \eqref{uno-cociclo4}
correspond bijectively to group homomorphisms $z: F/[F, F] \otimes
G/[G, G] \to Z(C)$, and the maps  $\chi: G \times F \to \widehat
C$ satisfying \eqref{uno-cociclo1} and \eqref{uno-cociclo2}
correspond bijectively to group homomorphisms $z: F/[F, F] \otimes
G/[G, G] \to \widehat C$.

As in Remark \ref{naif}, we look for $\sigma \in \Hom(G, Z^2(F,
\ku^{\times}))$ and $\tau: F \to Z^2(G, \ku^{\times})$, such that
$\tau_1 = 1$ and \begin{equation}\label{tau-naif}\tau_{xy}(t, s) =
\la \chi(t, g), z(s, x)\ra \, \tau_x(t, s) \, \tau_y(t, s),
\end{equation}
for all $x, y \in F$, $s, t \in G$.

\bigbreak A fairly complete answer can be given in the case when
$F$ and $G$ are cyclic. So assume that $F = \la a \ra$ has order
$N$ and $G = \la b \ra$ has order $M$. Let $C = \la u \ra$ be of
order $(M, N) = \gcd (M, N)$. Let also $\zeta \in \ku^{\times}$ be
a primitive $(M, N)$-th. root of unity. Define $z: G \otimes F \to
C$ and $\chi: G \otimes F \to \widehat C$ by
\begin{equation*}z(b^h \otimes a^j)  = u^{hj}, \qquad
\la \chi(b^h \otimes a^j), u^l \ra  = \zeta^{hjl};
\end{equation*} for all $0 \leq h \leq M-1$, $0 \leq j \leq N-1$,
$0 \leq l \leq (N, M)-1$.

\bigbreak We first determine the possible $\sigma$'s. Since $F$ is
cyclic, $Z^2(F, \ku^{\times}) = B^2(F, \ku^{\times})$; hence
giving $\sigma \in \Hom(G, Z^2(F, \ku^{\times}))$ is equivalent to
choosing $\sigma_b \in B^2(F, \ku^{\times})$ such that $\sigma_b^M
= 1$.

\bigbreak Concrete examples can be given as follows. Let $\omega
\in \ku^{\times}$ be such that $\omega^{MN} = 1$. We define $f: F
\to \ku^{\times}$ in the form $f(a^h) = \omega^h$, and let
$\sigma_b =
\partial f$ be the coboundary of $f$. So that
\begin{equation}\sigma_b(a^j, a^h) = \omega^{Nq},
\end{equation}
for all $0 \leq j, h \leq N-1$, where $j+h = Nq + r$, $0 \leq r
\leq N-1$.

\bigbreak We next consider the possibilities for $\tau$'s.

\begin{lema} The following are equivalent:

(i) $\tau: F \to Z^2(G, \ku^{\times})$ satisfies $\tau_1 = 1$ and
\eqref{tau-naif}.

(ii) There exists $\tau_a \in Z^2(G, \ku^{\times})$ satisfying
\begin{equation}\label{tau^N}\zeta^{st\frac{N(N-1)}{2}} \, \tau_a(b^s, b^t)^N
= 1,\end{equation} for all $0 \leq s, t \leq M-1$, such that
\begin{equation}\label{desc-tau}\tau_{a^m}(b^s, b^t) = \zeta^{st\frac{m(m-1)}{2}} \, \tau_a(b^s,
b^t)^m, \qquad 0 \leq s, t \leq M-1, \quad 0 \leq m \leq N-1.
\end{equation} \end{lema}

\pf (i) $\implies$ (ii). Condition \eqref{tau^N} follows from
$\tau_1 = 1$, and condition \eqref{desc-tau} follows from
\eqref{tau-naif} by induction on $m$.

(ii) $\implies$ (i). Left to the reader. \epf

\begin{obs} Let $\xi$ be a square root of $\zeta$, and let $\eta: G
\to \ku^{\times}$ be given by $\eta(b^s) = \xi^{s^2
\frac{N(N-1)}{2}}$. Condition \eqref{tau^N} can be rephrased as
saying that $\tau_a^N = \partial \eta$ in $Z^2(G, \ku^{\times})$.

If $(M, N)$ divides $\frac{N (N-1)}{2}$, then condition
\eqref{tau^N} amounts to $\tau_a^N = 1$.
\end{obs}

Let us give some concrete examples.

\begin{exa}\label{concrete1} Let $N, M \geq 1$ be integers such
that $(M, N)$ divides $\frac{N (N-1)}{2}$. Let $\zeta \in
\ku^{\times}$ be a primitive $(M, N)$-th. root of unity, and let
$\omega, \mu \in \ku^{\times}$ be such that $\omega^{MN} =
\mu^{MN} = 1$. Then there exists a braided Hopf algebra $R = \ku^G
{}^{\tau}\#_{\sigma}\ku F$, $F$ and $G$ as above, where $\fiz$,
$\fde$ are trivial and
\begin{align*}\sigma_{b^s}(a^j, a^h) & = \omega^{Nqs}, \quad
\text{if} \quad j+h = Nq + r, \, 0 \leq r \leq N-1,\\
\tau_{a^m}(b^s, b^t) & = \zeta^{st \frac{m(m-1)}{2}} \, \mu^{M
\tilde q m}, \quad \text{if} \quad s+t = M \tilde q + \tilde r, \,
0 \leq \tilde r \leq M-1.
\end{align*}
These braided Hopf algebras are commutative and cocommutative.

Note that we have $\sigma_{b^s}(a^j, a^h) = \nu(b^s, a^j) \nu(b^s,
a^h) \nu(b^s, a^{j+h})^{-1}$, for all $0 \leq s \leq M-1$, $0 \leq
j, h \leq N-1$, where $\nu(b^s, a^j) = \omega^{sj}$. Hence, the
braided Hopf algebra corresponding to the pair $(\omega, \mu)$ is
isomorphic to the braided Hopf algebra corresponding to the pair
$(1, \omega \mu)$. See Proposition \ref{prop-equiv}.
\end{exa}

\bigbreak The Examples given by Kashina in \cite{Ka} fit into the
present construction. Indeed, take $M = 2$, $N = 2^n$, with $n >
1$, $\zeta = -1$ and $\mu = 1$.

\bigbreak If $\omega = 1$, we get the braided Hopf algebra
$R^+_{n+1}$ in \cite{Ka}. This is dual to one of the examples in
\ref{nicosomm}.

\bigbreak If $\omega^N = -1$, we get the braided Hopf algebra
$R^-_{n+1}$ in \cite{Ka}. Again, this is dual to one of the
examples in \ref{nicosomm}.

\subsection{Cyclic groups.}

We shall  now consider the case where $F$ is a cyclic group of
order $N$. Write $F = \langle a: \, a^N = 1\rangle$. Let $f: G \to
\mathbb S_{N-1}$ be the group homomorphism associated to the
action $\fde : G \times F \to F$. By abuse of notation, we shall
use the same symbol to indicate an element $g \in G$ and its image
under $f$; so that we have $g \fde a^i = a^{g(i)}$, $g \in G$, $1
\leq i \leq N-1$.

\begin{lema}\label{diag-cyclic} Let $A$ be a finite abelian group.
The following collections of data are in bijective correspondence

(a) maps  $\alpha: G \times F \to A$ satisfying
\begin{align}
\label{uno-cociclo3-bis} \alpha(g, xy) & =  \alpha(g, x) \alpha(g
\fiz x, y), \qquad g\in G, x, y\in F, \\ \label{uno-cociclo4-bis}
\alpha(gh, x) & =  \alpha(g, h\fde x) \alpha(h, x), \qquad g, h\in
G, x\in F, \end{align} and

(b) maps  $\gamma: G \to A$ satisfying
\begin{align}
\label{gama1} 1 &  = \gamma(g) \gamma(g \fiz a) \dots \gamma(g
\fiz a^{N-1}), \\ \label{gama2} \gamma(gh) & = \gamma(g) \gamma(g
\fiz a) \dots \gamma(g \fiz a^{h(1)-1}) \gamma(h).
\end{align}
The bijection is given by
\begin{align}\label{c-gama}
\alpha(g, a^i) & =  \gamma(g) \gamma(g \fiz a) \dots \gamma(g \fiz
a^{i-1}), \qquad g \in G, \, 0 \leq i \leq N-1,\\ \alpha(1, a) & =
1. \end{align} \end{lema}

\pf Let $\gamma: G \to A$  be a function satisfying \eqref{gama1}
and \eqref{gama2}. Let $g \in G$, $0 \leq i, j \leq N-1$, and
write $i+j = Nq + r$, where $0 \leq r \leq N-1$ and $q = 0, 1$. By
definition, we have $$\alpha(g, a^{i+j}) = \gamma(g) \gamma(g \fiz
a) \dots \gamma(g \fiz a^{r-1}),$$ and on the other hand,
\begin{align*}
\alpha(g, a^i) \alpha(g \fiz a^i, a^j) &  = \gamma(g) \gamma(g
\fiz a) \dots \gamma(g \fiz a^{i-1}) \, \gamma(g\fiz a^i) \gamma(g
\fiz a^{i+1}) \dots \gamma(g \fiz a^{i+j-1}) \\ & = \gamma(g\fiz
a^N) \gamma(g \fiz a^{N+1}) \dots \gamma(g \fiz a^{N+r-1});
\end{align*}
the second equality because of  condition \eqref{gama1}. This
shows that $\alpha$ verifies \eqref{uno-cociclo3-bis}.

Let now $g, h \in G$, $0 \leq i \leq N-1$.  Observe that $gh \fiz
a^s = (g \fiz (h \fde a^s)) (h \fiz a^s) = (g \fiz a^{h(s)}) (h
\fiz a^s)$. Whence, in view of \eqref{gama2},
\begin{equation*}
\gamma(gh \fiz a^s) = \gamma((g \fiz a^{h(s)})  (h \fiz a^s)) =
\gamma(g \fiz a^{h(s)}) \dots \gamma(g \fiz a^{h(s)+ (h \fiz
a^s)(1)-1}) \gamma(h \fiz a^s).
\end{equation*}
Thus,
\begin{align*}
\alpha(gh, a^i) & = \gamma(gh) \gamma(gh \fiz a) \dots  \gamma(gh
\fiz a^{i-1}) \\ & = \gamma(g) \gamma(g \fiz a) \dots \gamma(g
\fiz a^{h(1)-1}) \gamma(h) \, \gamma(g \fiz a^{h(1)}) \dots
\gamma(g \fiz a^{h(1)+(h \fiz a)(1)-1}) \gamma(h \fiz a) \dots \\
& \gamma(g \fiz a^{h(i-1)}) \dots \gamma(g \fiz a^{h(i-1)+(h \fiz
a^{i-1})(1)-1}) \gamma(h \fiz a^{i-1}) \\ & =  \gamma(g) \dots
\gamma(g \fiz a^{h(i)-1}) \gamma(h) \dots \gamma(h \fiz a^{i-1}) =
z(h, h \fde a^i) z(h, a^i);
\end{align*}
the first equality by \eqref{gama2} and the third  equality
because of the following claim:

\begin{claim} We have $(h \fiz a^i) (1) = h(i+j)-h(i)$,
for all $i$. \end{claim}

\pf The  compatibility condition \eqref{comp1} implies that $(h
\fiz a^i) \fde a^j = (h \fiz a^i)^{-1}(h \fiz a^{i+j}) =
a^{h(i+j)-h(i)}$, whence the claim follows. \epf

Therefore, $\alpha$ satisfies \eqref{uno-cociclo4-bis}.

Conversely, assume   that the map $\alpha$ given by \eqref{c-gama}
satisfies \eqref{uno-cociclo3-bis} and \eqref{uno-cociclo4-bis}.
Using \eqref{uno-cociclo3-bis} and $a^N = 1$ to compute $\alpha(g,
a^N)$, one sees that $\gamma$ satisfies \eqref{gama1}. Similarly,
putting $x = a$ in \eqref{uno-cociclo4-bis}, the relation
\eqref{gama2} follows. This finishes the proof of the lemma. \epf

\begin{prop}  Let $\gamma: G \to \widehat C$ and $\eta: G \to Z(C)$ be maps satisfying \eqref{gama1} and \eqref{gama2}.
Let $\tau: G \times G \to (\ku^{\times})^F$ be a normalized 2-cocycle satisfying
\begin{equation}\tau_{a^{i+j}} (t, s) =
\la  \gamma(t \fiz a^{s(i)}) \, \gamma(t \fiz a^{s(i)+1}) \dots \gamma(t \fiz a^{s(i+j) - 1}), \eta(s) \eta(s \fiz a) \dots \eta(s \fiz a^{i-1}) \ra \,
\tau_{a^{i}} (t, s) \, \tau_{a^{j}} (t \fiz a^{s(i)}, s \fiz a^i),
\end{equation}
and $\tau_1(s, t) = 1$, for all $s, t \in G$, $0 \leq i, j \leq N-1$,
Then the bicrossed product $R = \ku^G {}^{\tau}\# \ku F$ is a braided Hopf algebra over $\ku C$, with respect to the maps $\chi: G \times F \to \widehat C$ and $z: G \times F \to Z(C)$ given by
\begin{align*}\chi(g, a^i) & = \gamma(g) \gamma(g \fiz a) \dots \gamma(g \fiz
a^{i-1}), \qquad \chi(1, a)  =
1, \qquad g \in G, \, 0 \leq i \leq N-1, \\
z(g, a^i) & = \eta(g) \eta(g \fiz a) \dots \eta(g \fiz
a^{i-1}), \qquad z(1, a)  =
1, \qquad g \in G, \, 0 \leq i \leq N-1.
\end{align*}

Moreover, all braided Hopf algebras admitting a diagonal realization over $\ku C$ are of this form. \end{prop}

\pf The first statement follows from Theorem \ref{braid-c-chi} and Lemma \ref{diag-cyclic}. See Remark \ref{naif}.
It follows also from Remark \ref{naif} that every braided Hopf algebra admitting a diagonal realization over $\ku C$ has this form. \epf

\section{Appendix}

The contents of this appendix have been suggested by the referee. It presents an alternative language which seems appropriate in discussions about matched pairs of groups.  This can be found for instance in \cite{maj-lib}; see also \cite{tak-ess-lyz}. The main constructions of the paper  are translated into this language. 

\subsection{Notation}
Let $\Sg$ be the set of all diagrams
$$
\begin{matrix} \quad g \quad \\ v \,\, \boxe \,\, x \\ \quad t \quad
\end{matrix},  \text{ where } g,t \in G, \, x, v \in F \text{ are such that } gx = vt.
$$ Thus $v = g \fde x$ and $t = g\fiz x$. Sometimes, we shall
simply write $\begin{matrix}  g \quad \\  \boxe \,\, x \\ \quad
\end{matrix} = \begin{matrix} \quad g \quad \\ v \,\, \boxe \,\, x \\ \quad t \quad
\end{matrix}$. A \emph{horizontal identity} is an element of the form
$\begin{matrix} 1 \quad \\  \boxe \,\, x \\ \quad
\end{matrix}$; a \emph{vertical identity} is an element of the form
$\begin{matrix} g \quad \\  \boxe \,\, 1 \\ \quad
\end{matrix}$.

Let $A = \begin{matrix} \quad g \quad \\ v \,\, \boxe \,\, x \\
\quad t \quad
\end{matrix}$, $B = \begin{matrix} \quad h \quad \\ w \,\, \boxe \,\, y \\ \quad s \quad
\end{matrix}$ be in $\Sg$. We shall write

\begin{flalign} & A \vert B \text{ if } x=w. \text{ Then } AB : =
\begin{matrix} \quad gh \quad \\ v \,\, \boxee \,\, y \\ \quad ts \quad
\end{matrix}  \text{ is in  } \Sg \text{ (horizontal product). }& \\
&\frac{A}{B} \quad \text{ if } t=h. \text{ Then } \begin{matrix}A
\\B\end{matrix} : =
\begin{matrix} \quad g \quad \\ vw \,\, \boxee \,\, xy \\ \quad s \quad
\end{matrix}  \text{ is in  } \Sg \text{ (vertical product). }&
\end{flalign}

\begin{obs} The notation $\begin{tabular}{p{0,4cm}|p{0,4cm}} $A$ & $B$ \\
\hline $C$ & $D$ \end{tabular}$ means that all possible horizontal
and vertical products are allowed; this implies that
$\displaystyle \frac{AB}{CD}$, $\begin{tabular}{p{0,4cm}|p{0,4cm}}
$A$ & $B$ \\ $C$ & $D$ \end{tabular}$ and there is no ambiguity in
the expression $\begin{matrix} A B \\ C D
\end{matrix}$. \end{obs}

\subsection{Cocycles} Let $\sigma$ and $\tau$ be as in \ref{bicpro}.
We define a function, that we still denote $\sigma$,  on the set
of all pairs $(A, B)$ with $\displaystyle \frac{A}{B}$,  and a
function $\tau$ on the set of all pairs $(A, B)$ with $A \vert B$,
by means of the formulas:

\begin{align}\sigma(A, B) & = \sigma_g(x, y), \\ \tau(A, B) & =
\tau_x(g, h), \end{align} where $A = \begin{matrix} g \quad
\\  \boxe \,\, x \\ \quad
\end{matrix}$ and $B = \begin{matrix}h \quad \\ \boxe \,\, y \\  \quad
\end{matrix}$ are in $\Sg$.
The cocycle and normalization conditions \eqref{cociclo-sigma},
\eqref{norm-sigma}, \eqref{cociclo-tau} and \eqref{norm-tau}
translate, respectively, as follows:
\begin{flalign} \label{cociclo-sigma-bis} & \text{If } \displaystyle
\frac{\displaystyle\frac{A}{B}}{C}, \text{ then } \sigma(A, B)
\sigma(\begin{matrix}A \\B\end{matrix},C) = \sigma(B, C) \sigma(A,
\begin{matrix}B \\C\end{matrix}). & \\
\label{norm-sigma-bis} & \text{If } A \text{ or } B \text{ is a
vertical identity, then } \sigma(A, B) = 1.& \\
\label{cociclo-tau-bis} & \text{If } A\vert B\vert C, \text{ then
} \tau(A, B) \tau(AB, C) =\tau(B, C) \tau(A, BC).& \\
\label{norm-tau-bis} & \text{If } A \text{ or } B \text{ is a
horizontal identity, then } \tau(A, B) = 1.&
\end{flalign}

\subsection{Operations} The bicrossed product $R = \ku^G \otimes \ku F$
has $\Sg$ as a basis with identification $\delta_g x =
\begin{matrix} g \quad \\   \,\, \boxe \,\, x \\ \quad \quad
\end{matrix}$. In this basis the operations of $R$ are determined by the
formulas

\begin{itemize}\item $A.B = \tau(A, B)
\begin{matrix}A
\\B\end{matrix}$,  if $\displaystyle \frac{A}{B}$,  and $0$ otherwise.

\item $\Delta(A) = \sum \sigma (B, C) B\otimes C$ , where the sum is
over all pairs $(B, C)$ with $B\vert C$ and $A = BC$.
\end{itemize}

\bigbreak The unitary conditions \eqref{norm2-sigma} and
\eqref{norm2-tau} translate into
\begin{flalign}
\label{norm2-sigma-bis} & \text{If } A \text{ or } B \text{ is a
horizontal identity, then } \sigma(A, B) = 1.& \\
\label{norm2-tau-bis} & \text{If } A \text{ or } B \text{ is a
vertical identity, then } \tau(A, B) = 1.&
\end{flalign}

\bigbreak When $A = \begin{matrix} \quad g \quad \\ x \,\, \boxe
\,\, y
\\ \quad h \quad
\end{matrix}$ is in $\Sg$, we put
\begin{align*} A^h &=  \begin{matrix} \quad g^{-1} \quad \\ y \,\, \boxd \,\, x \\
\quad h^{-1} \quad
\end{matrix} \quad\text{ (horizontal inverse), } \\
A^v &=  \begin{matrix} \quad h \quad \\ x^{-1} \,\, \boxe \,\, y^{-1} \\ \quad g \quad
\end{matrix} \quad\text{ (vertical inverse), } \\
A^{-1} &= (A^h)^v= (A^v)^h.
\end{align*}

The formula \eqref{antipode} for the antipode is then $$\Ss(A) =
\sigma(A^{-1}, A^h)^{-1} \tau(A^h, A)^{-1} A^{-1}.$$

\begin{obs} Note that $\begin{tabular}{p{0,7cm}|p{0,5cm}}
$A^{-1}$ & $A^v$ \\ \hline $A^h$ & $A$ \end{tabular}$ and
$\begin{tabular}{p{0,5cm}|p{0,7cm}}
 $A$ & $A^h$ \\ \hline  $A^v$ & $A^{-1}$ \end{tabular}$.
We have also $\sigma(A, A^v) = \sigma(A^v, A)$ from
\eqref{cociclo-sigma} and \eqref{norm-sigma}, and $\tau(A, A^h) =
\tau(A^h, A)$ similarly. \end{obs}

\subsection{Braiding} The braiding \eqref{def-c}  takes the
form $c(A \otimes B) = Q_{A, B} B \otimes A$, where $Q_{A, B} : =
Q^{x, y}_{g, h}$, for $A = \begin{matrix} g \quad \\  \boxe \,\, x
\\ \quad \end{matrix}$ and $B = \begin{matrix}h \quad \\ \boxe \,\, y \\  \quad
\end{matrix}$ in $\Sg$.

The compatibility condition \eqref{comp3} and the normalization
condition \eqref{norm-c} read, respectively, as follows:

\begin{flalign} \label{taktres} & \text{If } \begin{tabular}{p{0,4cm}|p{0,4cm}} A & B
\\ \hline C & D \end{tabular},  \text{ then  } \sigma(AB,
CD) \tau \left(\begin{matrix} A \\ C \end{matrix}, \begin{matrix}
B \\ D
\end{matrix}\right) = Q_{B, C} \, \tau(A, B) \tau(C, D) \sigma(A,
C) \sigma(B, D). & \\ &\text{If } A \text{ or } B \text{ is a
horizontal or vertical identity, then } Q_{A, B} = 1.&
\end{flalign}

Proposition \ref{const-struct} can be now stated as follows:

\begin{prop} For any $B$, $C$ in $\Sg$, there are unique $A$, $D$ in $\Sg$
with $\begin{tabular}{p{0,4cm}|p{0,4cm}} $A$ & $B$ \\ \hline $C$ &
$D$ \end{tabular}$. If we define $Q_{B, C}$ by \eqref{taktres},
then $R$ becomes a pre-braided Hopf algebra. \qed
\end{prop}

\subsection{Braided Hopf algebra} The compatibility conditions \eqref{comm1-m},
\eqref{comm2-m}, \eqref{comm1-delta}, \eqref{comm2-delta} in Lemma
\ref{1-cociclos} translate, respectively, as follows:
\begin{flalign} & \text{If } \displaystyle\frac{B}{C}, \text{ then } Q_{\displaystyle A,
\begin{matrix} B
\\ C\end{matrix}} = Q_{A, B}Q_{A, C}. & \\
&\text{If } \displaystyle\frac{A}{D}, \text{ then }
Q_{\displaystyle \begin{matrix} A
\\ D\end{matrix}, B} = Q_{A, B}Q_{D, B}. & \\
&\text{If } B \vert C, \text{ then } Q_{A, BC} = Q_{A, B}Q_{A, C}.
& \\ &\text{If } D \vert A, \text{ then } Q_{DA, B} = Q_{D,
B}Q_{A, B}. & \end{flalign}

\subsection{Realization} In Section \ref{diag-real}, the maps $z: \Sg \to Z(C)$ and $\chi: \Sg
\to \widehat C$ should satisfy the following conditions:
\begin{flalign}& \text{If } \displaystyle\frac{A}{B}, \text{ then } \chi\left(\begin{matrix}
A
\\ B\end{matrix}\right) = \chi(A)\chi(B). & \\
& \text{If } A \vert B, \text{ then } \chi(AB) = \chi(A)\chi(B).&
\\ & \text{If } \displaystyle\frac{A}{B}, \text{ then }
z\left(\begin{matrix} A \\ B\end{matrix}\right) = z(A)z(B). & \\ &
\text{If } A \vert B, \text{ then } z(AB) = z(A)z(B). &
\end{flalign}

Finally, equation \eqref{tak-comp} is now $Q_{A, B} = \langle
\chi(B), z(A)\rangle$.


\begin{thebibliography}{999}

\bibitem[EG]{eg-exp}  P. Etingof and S. Gelaki,
\emph{On the exponent of finite dimensional Hopf algebras},
 Math. Res. Lett. \textbf{6}   (1999), pp. 131--140.

\bibitem[K]{kac}  G. I. Kac,
\emph{Extensions of groups to ring groups}, Math. USSR Sbornik
{\bf 5} (1968), pp.  451--474.

\bibitem[Ka]{Ka}  Y. Kashina,
\emph{Examples of  Hopf algebras of dimension $2^m$},
MSRI Preprint  \#2000-003.

\bibitem[Mj]{maj} S. Majid,  \emph{Physics for algebraists:
Non-commutative and non-cocommutative Hopf algebras by a
bicrossproduct construction}, J.  Algebra {\bf 130} (1990), pp.
17--64.

\bibitem[Mj2]{maj-lib} S. Majid,  \emph{Foundations of quantum group theory}, Cambridge Univ. Press, Cambridge (1995).

\bibitem[M]{Maext}  A. Masuoka,
\emph{Extensions of  Hopf algebras},
 Trabajos de Matem\' atica \textbf{41/99}, Fa.M.A.F.   (1999).

 Available at {\tt http://www.mate.uncor.edu/andrus}.

\bibitem[N]{pqq} S. Natale,  \emph{ On semisimple Hopf algebras of dimension
$pq^2$},   J.  Algebra {\bf 221} (1999), pp.  242--278.

\bibitem[Sb]{Sb}  P. Schauenburg,
\emph{On the braiding on a Hopf algebra in a braided category},
 New York J. Math. \textbf{4}   (1998), pp. 259--263.

\bibitem[S]{S}  Y. Sommerh\"auser,
\emph{Yetter-Drinfel'd Hopf Algebras over Groups of Prime Order},
Lecture Notes in Math. \textbf{1789}   (2002), Springer-Verlag.

\bibitem[T1]{tak}  M. Takeuchi,
\emph{Matched pairs of groups and bismash products of Hopf
algebras}, Commun. Algebra \textbf{9}   (1981), pp. 841--882.

\bibitem[T2]{T}  M. Takeuchi,
\emph{Survey of braided Hopf algebras},
 Contemp. Math. \textbf{267}   (2000), pp. 301--324.
 
\bibitem[T3]{tak-ess-lyz}  M. Takeuchi,
\emph{Survey on matched pairs of groups. An elementary approach to the ESS-LYZ theory}, preprint (2001).

\end{thebibliography}
\end{document}